\documentclass[12pt,oneside]{amsart}
\usepackage{amssymb,amsmath}
\makeatletter
\def\@cite#1#2{{\m@th\upshape\bfseries%
[{#1\if@tempswa{\m@th\upshape\mdseries, #2}\fi}]}}
\makeatother
\theoremstyle{plain}
\newtheorem{thm}{Theorem}[section]
\newtheorem{lem}[thm]{Lemma}
\newtheorem{cor}[thm]{Corollary}
\newtheorem{prop}[thm]{Proposition}
\theoremstyle{definition}
\newtheorem{rem}[thm]{Remark}
\newtheorem{defn}[thm]{Definition}
\newtheorem{note}[thm]{Note}
\newtheorem{eg}[thm]{Example}

\newcommand{\Prf}{\noindent\textbf{Proof.\ }}
\newcommand{\bx}{\strut\hfill$\blacksquare$\medbreak}
\newcommand{\upbx}{\vspace{-2.5\baselineskip}
\newline\hbox{}\hfill$\blacksquare$\newline\medbreak}

\newcommand{\BH}{{\B(\H)}}
\newcommand{\ca}{\mathrm{C}^*}

\newenvironment{sbmatrix}{\left[\begin{smallmatrix}}{\end{smallmatrix}\right]}

\DeclareMathOperator*{\wotlim}{\textsc{wot}--lim}

\DeclareMathOperator*{\sotlim}{\textsc{sot}--lim}

\newcommand{\wot}{\textsc{wot}}

\newcommand{\bbC}{{\mathbb{C}}}

  \newcommand{\B}{{\mathcal{B}}}

  \newcommand{\E}{{\mathcal{E}}}
  \newcommand{\F}{{\mathcal{F}}}
  
\renewcommand{\H}{{\mathcal{H}}}

  \newcommand{\K}{{\mathcal{K}}}

  \newcommand{\M}{{\mathcal{M}}}
  
\renewcommand{\O}{{\mathcal{O}}}

  \newcommand{\V}{{\mathcal{V}}}
  \newcommand{\W}{{\mathcal{W}}}
  \newcommand{\X}{{\mathcal{X}}}

\newcommand{\eps}{\varepsilon}

\newcommand{\upchi}{{\raise.35ex\hbox{$\chi$}}}
\newcommand{\fL}{{\mathfrak{L}}}

\newcommand{\fR}{{\mathfrak{R}}}

\newcommand{\Lat}{\operatorname{Lat}}
\newcommand{\ran}{\operatorname{Ran}}
\newcommand{\rk}{\operatorname{rk}}
\newcommand{\spn}{\operatorname{span}}
\newcommand{\tr}{\operatorname{tr}}

\newcommand{\qand}{\quad\text{and}\quad}

\newcommand{\qfor}{\quad\text{for}\quad}

\newcommand{\crv}{\operatorname{K}}
\newcommand{\elr}{\operatorname{\upchi}}
\newcommand{\prk}{\operatorname{pure \, rank}}

\newcommand{\ktil}{\tilde{\crv}}

\title{\normalsize\bf \rule{0pt}{3.5cm} THE CURVATURE INVARIANT OF A \\
NON-COMMUTING $N$-TUPLE}

\author{DAVID W. KRIBS}
\thanks{The author was partially supported by an Ontario Graduate Scholarship}
\date{}
\begin{document}
\maketitle
\thispagestyle{empty}

\baselineskip=12pt
\begin{quote}
Non-commutative versions of Arveson's curvature invariant and Euler characteristic for a
commuting $n$-tuple of operators are introduced. The non-commutative curvature invariant is
sensitive enough to
determine if an $n$-tuple is free.  In general both invariants can be thought of as measuring the
freeness or curvature of an $n$-tuple. The connection with dilation theory provides motivation
and exhibits relationships between the invariants.
A new class of examples is used to illustrate the differences encountered in the
non-commutative setting and obtain information on the ranges of the invariants. The curvature
invariant is also shown to be upper semi-continuous.
\end{quote}

\vskip 1truecm
\baselineskip=15pt

A notion of  curvature was recently introduced in commutative multi-variable operator theory by
Arveson \cite{Arv_euler}. Specifically, the {\it curvature invariant} and {\it Euler
characteristic} for a contractive $n$-tuple of mutually commuting operators were developed
from two different perspectives, both dependent on commutative methods. In addition,
asymptotic
formulae were obtained for the invariants in terms of a certain completely positive map.
It turns out that both invariants are always integers. In fact, for the classes which Arveson
focuses on the two perspectives yield the
same invariant.
The goal of this paper is to develop non-commutative versions of these invariants.
The non-commutative versions introduced here can be thought of as the {\it analogues} of
Arveson's to the non-commutative setting. Indeed, this claim is well supported since they
possess many of the basic properties of their
commutative relatives. Most importantly, the basic property of the non-commutative curvature
invariant is that it is
sensitive enough to detect when an $n$-tuple is free.  However, there are definite differences in
the non-commutative setting. A new class of examples is
included which help illustrate these differences. For instance, there are many interesting cases
when the invariants are distinct and not integers. In fact the range of the curvature invariant
is the entire positive real line. Further, a complete continuity picture for the
invariants is
painted. Namely, the curvature invariant is shown to be upper semi-continuous with respect to
the natural notions of convergence.

The first section contains requisite preliminary material. This includes a discussion of the
completely positive map defined by every row contraction which determines the
invariants. The related dilation theory is also recalled. In the second section, the existence of the
non-commutative curvature invariant and Euler characteristic is proved. The connection with
dilation theory is used to provide motivation and establish a hierarchy of the invariants.

The third section contains an analysis of how the invariants behave on pure contractions.  The
basic
property is that the curvature invariant can detect when an $n$-tuple is free. In other words the
$n$-tuple is  unitarily
equivalent to copies
of the pure isometries. This is analogous to Arveson's basic property;
however, different methods must be used to prove the result here. Also, an analogous
characterization of when the curvature invariant annihilates a contraction is obtained. For pure
contractions, the Euler
characteristic can provide a measure of the freeness of an $n$-tuple.

In the fourth section,  continuity and stability properties are investigated. The curvature invariant
is upper semi-continuous; however, the rigidity of the Euler characteristic prevents any
non-trivial continuity results. The invariants are shown to be stable under certain conditions, but
not to the same extent as Arveson's.

In the last section, a new class of pure contractions is introduced which yield information on the
ranges of the invariants.  This class consists of  finite rank perturbations of a subclass of the
atomic
Cuntz representations originally considered by K. Davidson and D. Pitts in \cite{DP2}.
The curvature invariant and Euler characteristic are shown to be distinct and non-integers in
general. Two
different collections of examples are used to prove that the range of the curvature invariant is the
entire positive real line. Further, various open problems related to the ranges are pointed out.

The content of this paper is part of the doctoral thesis \cite{Kri_thesis}. The author would like to
thank his advisor Ken Davidson for his assistance.

\begin{note}
Upon submission of this paper the author became aware of the recently submitted paper
\cite{Pop_curv} by G. Popescu. There is a lot of overlap between our papers, although the
perspectives and proofs are quite different. In fact Popescu obtains the same curvature invariant.
Further, he obtains a pair of pseudo-Euler characteristics which both turn out to be equal to the
Euler characteristic developed here. This comes as a corollary of Theorem ~\ref{lim_exist}.
Following the non-commutative analogue of Arveson's algebraic approach used in
\cite{Arv_euler}, Popescu is also able to prove a non-commutative version of Arveson's
Gauss-Bonnet-Chern Theorem and show that the range of the Euler characteristic is the positive
real
line. Lastly, an elementary class of examples showing how to obtain the positive real line as the
range of both invariants, which the author discovered after submission, has been included in
Note ~\ref{cuteeg}.
\end{note}

%%%%%%%%%%%%%
\section{Preliminaries} \label{S:prelim}
%%%%%%%%%%%%%

To every row contraction $A = (A_1, \ldots, A_n)$ of operators $A_i$ acting on a Hilbert space
$\H$, there is a corresponding completely positive map $\Phi (\cdot)$ on $\BH$ defined by
\[
\Phi(X) = \sum_{i=1}^n A_i X A_i^* = AX^{(n)}A^*.
\]
When there is more than one $n$-tuple involved, the associated map will be denoted $\Phi_A
(\cdot)$.  This map is also completely contractive since
\[
\Phi(I) = \sum_{i=1}^n A_i  A_i^* = A A^* \leq I.
\]
Given a word $w = i_1 \cdots i_d$ in the unital free semigroup on $n$ generators $\F_n$, define
a contraction $A_w := w(A) = A_{i_1}\cdots A_{i_d}$. For $k \geq 1$, the
sequence
\[
\Phi^k (I) = \sum_{i=1}^n A_i \Phi^{k-1}(I) A_i^* = \sum_{|w| =k} A_w A_w^*
\]
satisfies $\Phi(I) \geq \Phi^2 (I) \geq \ldots \geq 0$. It is this sequence of operators which
Arveson uses to derive information on the original $n$-tuple. As a decreasing sequence of
positive contractions, the sequence has a strong operator
topology limit which has been denoted by
\[
\Phi^\infty (I) = \sotlim_{k\rightarrow\infty}  \Phi^k (I).
\]
The two extreme cases are important.

\begin{defn}
Let $A = (A_1, \ldots, A_n)$ be a contractive $n$-tuple of operators acting on $\H$. Then
\begin{itemize}
\item[$(i)$] $A$ is a {\bf pure} contraction if $\Phi^\infty (I) = 0$
\item[and]
\item[$(ii)$] $A$ is a {\bf Cuntz} contraction if $\Phi^\infty (I) = I$.
\end{itemize}
\end{defn}

These definitions are motivated by the special case of an $n$-tuple of isometries $S = (S_1,
\ldots, S_n)$. In fact much
of the theory here is motivated by this special case. The isometries $S_i$
have pairwise orthogonal ranges exactly when the $n$-tuple is contractive. That is,
\[
\sum_{i=1}^n S_iS_i^* \leq I \quad\text{if and only if}\quad S_i^*S_j = \delta_{ij}I.
\]
The important
example in the pure case, and in this paper,  is determined by the left regular representation
$\lambda$ of $\F_n$ on $n$-variable
Fock
space, which can be viewed as  $\H_n := \ell^2 (\F_n)$. From this perspective $\H_n$ has an
orthonormal basis given by $\{ \xi_w: w\in\F_n \}$. The
isometries $L_i = \lambda (i)$ are defined by $L_i \xi_w = \xi_{iw}$, and also arise in
physics as the
`left creation operators'. Similarly, isometries
determined by words $w$ are denoted  $\lambda(w) = L_w$. The $n$-tuple $L= (L_1, \ldots,
L_n)$ satisfies
\[
L L^* = \sum_{i=1}^n L_i L_i^* = I - P_e,
\]
where $P_e$ is the orthogonal projection onto the span of the so called vacuum vector $\xi_e$
determined by the identity $e$ of $\F_n$ (which
corresponds to the empty word). The defect projections associated with $L$ are key in the
analysis.
\begin{defn}
For $k \geq 1$, let
\[
Q_k := I - \Phi_L^k (I) = I - \sum_{|w|=k} L_w L_w^*
\]
be the orthogonal projection onto $\spn \{ \xi_w : |w| < k \}$.
\end{defn}

There has been extensive recent work related to these operators. In
particular, the $\wot$-closed non-selfadjoint algebras $\fL_n$ generated by the $L_i$
have been shown by Davidson, Pitts and Arias, Popescu to be the appropriate non-commutative
analytic Toeplitz algebras
\cite{AP,DP1,DP2,Pop_diln,Pop_beur}. Also see \cite{AP2,DP3,Kri,Pop_fact} for more
detailed information on these
algebras. The \wot-closed algebras corresponding to the right regular representation are denoted
by $\fR_n$.
The generating isometries are defined by $\rho(w)=R_{w'}$ where $R_u \xi_v=\xi_{vu}$, and
$w'$ denotes the word $w$ in reverse order. It is unitarily equivalent to $\fL_n$ and is precisely
the commutant of $\fL_n$.

There is a nice link between this example and Arveson's motivating example. The key
contraction
in the commutative setting is given by multiplication $M_i$ for $1 \leq i \leq n$, by the $n$
coordinate functions on $n$-variable symmetric Fock space $\H_n^s$. In terms of the notation
here, this space can be realized as the subspace of $\H_n$ spanned by unit vectors of the form
$\frac{1}{k!} \sum_{\sigma\in S_k} \xi_{\sigma (w)}$ for $w$ in $\F_n$ with $|w| = k$
($\sigma$ is the permutation in the symmetric group $S_k$ which permutes the terms of $w$).
Alternatively, in
\cite{DP2} the complement $(\H_n^s)^\perp$ was shown to be the full range in $\Lat \fL_n \cap
\Lat \fR_n$ of the $\wot$-closed commutator ideal in $\fL_n$. It is given by
\[
(\H_n^s)^\perp = \spn \{ \xi_{uijv} - \xi_{ujiv}: u,v\in\F_n \}.
\]
Since $M = (M_1, \ldots, M_n)$ is a pure contraction with $\rk (I - \Phi_M (I))=1$, it follows
that $L$ is the minimal isometric dilation of $M$ (see below). Thus, $\H_n^s$ is co-invariant
for the $L_i$ and each $M_i = P_{\H_n^s} L_i |_{\H_n^s} = (L_i^*|_{\H_n^s})^*$.
This can be proved directly, but it is automatic from the dilation theory discussed below.

For an $n$-tuple of isometries $S = (S_1, \ldots, S_n)$ with pairwise orthogonal ranges, the
subspaces
\[
\H_p := \{ x\in\H : \lim_{k\rightarrow\infty} (\Phi^k(I) x, x) = 0 \}
\]
and
\[
\H_c := \{ x\in\H : \Phi^k(I) x = x , \,\,\,{\rm for}\,\,\, k \geq 1 \}
\]
can be seen to reduce the $S_i$ and are orthogonal complements of each other. This fact
leads to Popescu's
Wold decomposition \cite{Pop_diln}. The  $S_i$ simultaneously decompose as a direct sum
$S_i \simeq T_i \oplus L_i^{(\alpha)}$, where the $T_i = S_i|_{\H_c}$ are isometries which
form a representation of the Cuntz C$^*$-algebra $\O_n$ since $\sum_{i=1}^n T_i T_i^* =
I_{\H_c}$. The
multiplicity of the left regular representation is given by the rank of the projection $I -
\sum_{i=1}^n S_i S_i^*$. This follows because the range $\W$ of this projection is {\it
wandering} for $S$. That is, the subspaces $S_v\W$ are pairwise orthogonal for distinct words
$v$ in $\F_n$. Although Cuntz and pure contractions are important, there is no such
decomposition for an arbitrary contraction. The point is that for a general contractive $n$-tuple
$A$, the subspaces $\H_p$ and $\H_c$ will not be orthogonal complements as they can
be skewed when the $A_i$ are not isometries. An example of this phenomenon was provided in
\cite{DKS}.

The connection with dilation theory is also important. Analogous to the unique minimal
isometric
dilation of a contraction to an isometry established by Sz.-Nagy \cite{SF}, every contractive
$n$-tuple of operators has recently been shown to have a unique minimal joint isometric dilation
to an
$n$-tuple of isometries with pairwise orthogonal ranges acting on a larger space. This  theorem
derives from the work of
Frahzo, Bunce and Popescu \cite{Fra1,Bun,Pop_diln}. It is this dilation which will be considered
throughout the paper. Let $A = (A_1, \ldots, A_n)$ be a
contraction acting on $\H$, with minimal isometric dilation $S = (S_1, \ldots, S_n)$ acting
on $\K$. Then each $A_i$ is the compression to $\H$ of $S_i$. Minimality means the closed
span of the subspaces $S_w \H$ is all of $\K$. The condition of uniqueness is up to unitary
equivalence fixing $\H$. Further, from the construction of the dilation, the subspace $\H$ of
$\K$ is co-invariant for
$S$. In fact, the following spatial decomposition of $\H^\perp$ was obtained in \cite{DKS} as
Lemma 3.1.

\begin{lem}\label{wander}
The subspace $\W = (\H + \sum_{i=1}^n S_i \H ) \ominus \H$ is a wandering subspace for $S$,
and $\sum_{v\in\F_n} \oplus S_v \W = \H^\perp$.
\end{lem}

Thus, $\H^\perp$ is unitarily equivalent to a multiple $\H_n^{(\alpha)}$ of Fock space, where
$\alpha = \dim \W$, and $S_i|_{\H^\perp} \simeq L_i^{(\alpha)}$.  Hence decomposing $\K =
\H \oplus \H^\perp$, each $S_i$ may be written as a matrix $S_i \simeq
 \begin{sbmatrix} A_i & 0 \\ X_i & L_i^{(\alpha)}\end{sbmatrix}$.

Popescu's Wold decomposition shows that every contractive $n$-tuple $A$ determines a Cuntz
representation and copies of the left regular representation through its minimal isometric dilation
$S$. Every Cuntz contraction has a Cuntz representation as its dilation. For, it is
not hard to see that $\sum_{i=1}^n A_i A_i^* =I$ if and only if $\sum_{i=1}^n S_i S_i^* =I$.
When an $n$-tuple $A$ acts on finite dimensional space, the associated Cuntz representations
from the dilations
have been completely classified \cite{DKS}. As an aside, overall there has been considerable
recent interest  in classifying Cuntz representations.
For instance, see \cite{Arv_III,BJP,BJendII,BJit,DKS,Kri_cpdiln,Laca,Pow}. In addition to the
connection
with dilation theory, there is also a
 correspondence between Cuntz representations and endomorphisms of $\BH$. Further, certain
Cuntz representations give rise to wavelets \cite{BJwave}. One of the examples in
Section~\ref{S:examples} is
related to such a representation.

However, the strength of this paper is in the information it yields for pure contractions. Recall
that the multiplicity of
the left regular representation in the dilation $S$ of $A$ is given by the rank of $I -
\sum_{i=1}^n S_i S_i^*$.

\begin{defn}
Let $A = (A_1, \ldots, A_n)$ be a contraction with minimal isometric dilation $S = (S_1, \ldots,
S_n)$. Then define the {\bf pure rank} of $A$ to be the multiplicity
\[
\prk(A) = \rk \Bigl( I - \sum_{i=1}^n S_i S_i^* \Bigr).
\]
\end{defn}

The pure rank is a key invariant in this paper and can be computed directly in terms of the
original
$n$-tuple. It was introduced and described in \cite{DKS} for finite dimensional $n$-tuples. The
proof actually
works for an arbitrary contraction and is included here for the sake of completeness.

\begin{lem}\label{pure_rank}
The pure rank of a contractive $n$-tuple $A$ is computed as
\[
  \prk( A )
              = \rk \Bigl( I - \sum_{i=1}^n A_iA_i^* \Bigr) = \rk( I - \Phi(I) ).
\]
\end{lem}

\Prf  The wandering space is
$\X = \ran \bigl( I -  \sum_{i=1}^n S_iS_i^* \bigr)$
and the pure rank of $A$ equals $\dim \X$.
The minimality of the dilation means that $\X$ does not intersect
$\H^\perp$.
Therefore $P_\H P_\X P_\H$ has the same rank as $P_\X$.
However, observe that
\begin{align*}
  P_\H P_\X P_\H |_\H
  &= P_\H \Bigl( I_\K - \sum_{i=1}^n S_iS_i^* \Bigr) P_\H |_\H \\
  &=  I_\H - \sum_{i=1}^n A_iA_i^* = I_\H - \Phi(I_\H) .
\end{align*}
Thus $\prk(A) = \rk \bigl( I - \Phi(I) \bigr)$.
\bx

This quantity is what Arveson calls the {\it rank} of a contraction. The change in terminology
here to {\it pure rank} is motivated by the link with dilation theory.
A valuable property of pure contractions is that they are compressions of multiples of the left
regular representation, the pure isometries. This was first observed by Popescu \cite{Pop_diln}
and is included here
to illustrate the connection which the minimal dilation can have with the original $n$-tuple.

\begin{prop}
If $A = (A_1, \ldots, A_n)$ is a non-zero pure contraction, then the minimal isometric dilation of
$A$ is $L^{(\alpha)} := (L_1^{(\alpha)}, \ldots, L_n^{(\alpha)})$, where $\alpha = \prk(A)$.
\end{prop}

\Prf
Suppose $A$ acts on $\H$ and its minimal dilation $S$ acts on $\K$. It suffices to show that
$\K = \K_p$. If this is the case, then $S \simeq L^{(\alpha)}$ where $\alpha = \prk(A)$.

Consider a vector of the form $S_vx$ for some $x$ in $\H$ and $v$ in $\F_n$. Then for
$k > |v|$,
\begin{eqnarray*}
|| \Phi^k_S (I) S_v x ||^2 &=& \sum_{|w|=k} || S_w S_w^* S_v x ||^2 \\
          &=& \sum_{|u| = k - |v|} || S_u^* x ||^2 \\
          &=& (\Phi^{k - |v|}_A (I) x, x),
\end{eqnarray*}
which converges to zero since $A$ is pure. By minimality it follows that $\Phi^\infty_S(I) = 0$,
so that $\K = \K_p$.
\bx

In  \cite{Arv_euler}, Arveson considers contractive $n$-tuples for which the $A_i$ are pairwise
commuting and $I -
\Phi(I)$ is finite rank (ie: $\prk(A) < \infty$). The invariants introduced in this paper can be
thought of as the non-commutative {\it analogues} of Arveson's. Indeed, a strong case is made
for
this claim in Sections~\ref{S:curv_inv} and ~\ref{S:basic_prop}. They are defined
for any contraction; however, as in Arveson's setting the greatest amount of
information is obtained in the
finite rank case. In this case, the defect operators
\[
I - \Phi^k(I) = \sum_{i=0}^{k-1} \Phi^i (I - \Phi(I))
\]
form an increasing sequence of finite rank positive contractions. The idea is to use this sequence
to
obtain information on the associated $n$-tuple. Arveson defines two invariants from different
perspectives. The curvature invariant is defined by integration of the trace of a certain defect
operator which is defined on the range of $I - \Phi(I)$. On the other hand, if $A$ is a commuting
$n$-tuple acting
on $\H$, then a natural commutative Hilbert $A$-module structure can be placed on $\H$.  The
$A$-submodule determined by the range of $I - \Phi(I)$ is finitely generated, hence from
commutative module theory this module has a finite free resolution. The ranks of the free
modules from the associated exact sequence are known as the `Betti numbers'. The Euler
characteristic is defined as the alternating sum of these ranks. An operator theoretic version of
the Gauss-Bonnet-Chern theorem is obtained from this point of view.

The following asymptotic formulae were obtained for the commutative curvature invariant and
Euler characteristic respectively:
\[
n! \lim_{k\rightarrow\infty} \frac{\tr (I - \Phi^k (I))}{k^n}
\qand
n! \lim_{k\rightarrow\infty} \frac{\rk (I - \Phi^k (I))}{k^n}.
\]
Of course, when one moves to the non-commutative setting the strong links with commutative
module and function theory are lost. However, the map $\Phi_A$ is a bonafide completely
positive map for {\it any} $n$-tuple. Hence this seems like the natural point of view to take
when attempting to develop the non-commutative versions. To obtain these new
invariants, it becomes clear that the traces and ranks
of the operators $I - \Phi^k(I)$ must be re-normalized. Indeed, it is easy to obtain examples in
the general setting, even in the pure rank one case, for which Arveson's invariants are infinite.
So the factors $k^n / n!$ turn out to be
specific to the commutative setting.

Throughout the paper, $n$ is taken to be a finite integer with $n \geq 2$. Although the results
used  from dilation theory and the theory of non-commutative analytic Toeplitz algebras go
through for $n = \infty$, the invariants considered here are not defined in this case. However,
there may be analogous invariants in the $n = \infty$ setting.

%%%%%%%%%%%%%%%
\section{The Non-commutative Invariants} \label{S:curv_inv}
%%%%%%%%%%%%%%%

In this section, the non-commutative versions of the curvature invariant and Euler characteristic
are developed. The connection with dilation theory is utilized to provide motivation for the
definitions and leads to a general hierarchy of the related invariants. The starting point is an
elementary lemma.

\begin{lem}
If $\{ a_k \}_{k\geq 1}$ and $\{ s_k \}_{k\geq 1}$ are sequences of non-negative numbers with
$a_{k+1} \leq a_k + s_k$ and $\sum_{k\geq 1} s_k < \infty$, then $\lim_{k\rightarrow\infty}
a_k$ exists.
\end{lem}

\Prf
Let $\alpha_1 = \lim \inf a_k$ and $\alpha_2 = \lim \sup a_k$. Suppose $0\leq \alpha_1 <
\alpha_2$, with $\alpha_2 - \alpha_1 \geq \delta > 0$. Choose positive integers $m_1 >
n_1$ such that $a_{m_1} - a_{n_1} \geq \delta /2$. Then one has
\begin{eqnarray*}
\delta /2 &\leq&  a_{m_1}-a_{m_1-1}+a_{m_1-1}-\ldots +
a_{n_1+1}-a_{n_1} \\
                                    &\leq& \sum_{n_1\leq k \leq m_1-1} s_k.
\end{eqnarray*}
In a similar fashion obtain a sequence of positive integers $n_1<m_1<n_2<m_2<\ldots$ for
which $\sum_{n_j\leq k \leq m_j-1} s_k \geq \delta/2$ for $j\geq 1$. This would contradict the
summability of the sequence $s_k$.
\bx

Note that this lemma does not hold if the sequence $s_k$ simply converges to zero. With
this result in hand, the existence of the invariants in the general non-commutative setting can be
proved. At first glance the proof may seem somewhat sterile. This subtlety  is clarified by the
situation for an $n$-tuple of isometries
which is described below.

\begin{thm} \label{lim_exist}
Let $A=(A_1, \ldots ,A_n)$ be a contraction acting on $\H$. Then the
limits
\begin{itemize}
\item[$(i)$] $(n-1) \lim_{k\rightarrow\infty}\tr (I - \Phi^k (I)) / n^k$
\item[{\it and}]
\item[$(ii)$] $(n-1) \lim_{k\rightarrow\infty} \rk (I - \Phi^k (I)) / n^k,$
\end{itemize}
both exist.
\end{thm}

\Prf
The key identity used to prove the existence of both limits is
\begin{eqnarray}\label{key_id}
I - \Phi^{k+1}(I) &=& I - \Phi(I) + \Phi(I- \Phi^k(I)).
\end{eqnarray}
To prove $(i)$, notice that if $X$ is a positive operator on $\H$, then one has
\begin{eqnarray*}
\tr (\Phi(X)) &=& \sum_{i=1}^n \tr (A_i X A_i^*) \\
          &=& \sum_{i=1}^n \tr (X^{1 / 2}A_i^*A_i X^{1 / 2}) \leq n \tr(X).
\end{eqnarray*}
Observe that since the sequence $I - \Phi^k(I)$ is increasing, the limit of the normalized traces is
trivially infinite when $\tr (I - \Phi(I)) = \infty$. Otherwise, equation $(1)$ yields for $k \geq 1$,
\begin{eqnarray*}
0 \leq \frac{\tr (I-\Phi^{k+1}(I))}{n^{k+1}} \leq \frac{\tr (I-\Phi^{k}(I))}{n^{k}} + \frac{\tr (I -
\Phi(I))}{n^{k+1}}.
\end{eqnarray*}
Hence the lemma applies and the existence of the associated limit is proved.

The proof of $(ii)$ proceeds in a similar manner. If
\[
\rk (I - \Phi(I)) = \prk(A) = \infty,
\]
it follows that the limit of the ranks is trivially infinite. In the finite pure rank case, basic
properties of the rank function prove that
\begin{eqnarray*}
0 \leq \rk (I-\Phi^{k+1}(I)) &\leq& \rk (\Phi ( I-\Phi^{k}(I))) + \rk (I-\Phi (I)) \\
               &\leq& n \rk ( I-\Phi^{k}(I)) + \prk (A).
\end{eqnarray*}
Dividing by $n^{k + 1}$ shows that the lemma can be applied again.
\bx

\begin{rem}
These limits exist for any contractive $n$-tuple; the point is that they yield good information
when the pure rank, the rank of the defect operator $I - \Phi(I)$, is finite. These are the
contractions which will be focused on throughout the rest of the paper.
The naive motivation for considering these particular re-normalizations comes from an analysis
of the words in $n$ letters. Consider the words on $n$ letters of length less than $k$ for large
$k$. The total number of words in $n$ commuting letters is on the order of $k^n / n!$, whereas
the
number of words in $n$ non-commuting letters is on the order of $n^k / (n-1)$. This
is clarified further below by considering the situation for an $n$-tuple of isometries with
pairwise orthogonal ranges. Lastly, the uniformity of the convergence estimates obtained in the
theorem allow one to obtain continuity results. These are established in  Section
~\ref{S:cont_stab}.
\end{rem}

In keeping with Arveson's nomenclature, the notation  from the commutative setting will be
preserved.

\begin{defn}
Let $A=(A_1, \ldots, A_n)$ be a contraction.
\begin{itemize}
\item[$(i)$] The {\bf curvature invariant} of $A$ is defined to be the limit
\[
\crv (A) := (n-1) \lim_{k\rightarrow\infty} \frac{\tr (I-\Phi^k_A (I))}{n^{k}},
\]
\item[and]
\item[$(ii)$] The {\bf Euler characteristic} of $A$ is defined to be the limit
\[
\elr (A) := (n-1) \lim_{k\rightarrow\infty} \frac{\rk (I-\Phi^k_A(I))}{n^{k}}.
\]
\end{itemize}
\end{defn}

Thus these invariants are {\it not} generalizations of Arveson's; instead, they are the analogues
in the non-commutative setting. This point is supported by the properties developed throughout
the rest of the paper. An obvious property which is the same is the inequality
$\crv(A) \leq \elr(A)$. It holds for the same reason here; the operators $I - \Phi^k (I)$ are
positive contractions. As discussed in the preliminary section, $n$-tuples of isometries provide
much of the motivation.
Before considering them a simple but helpful lemma is presented.

\begin{lem}\label{dir_sum}
Let $A \oplus B$ be the coordinate-wise direct sum of two contractive $n$-tuples $A$ and $B$.
Then $\crv (A \oplus B) = \crv (A) + \crv (B)$ and
$\elr (A \oplus B) = \elr (A) + \elr (B)$.
\end{lem}

\Prf
This follows from the identity
\[
\tr (I - \Phi^k_{A \oplus B} (I)) = \tr  (I - \Phi^k_{A} (I)) + \tr (I - \Phi^k_{B} (I)).
\]
The same is true for the ranks.
\bx

%As an immediate corollary, it becomes concretely apparent that these invariants yield
%information on the pure part of a given $n$-tuple.
%
%\begin{cor}
%Let $A$ be a contractive $n$-tuple with its decomposition into Cuntz and pure parts $A = A_c
%\oplus A_p$. Then
%$\crv (A) = \crv (A_p)$ and $\elr (A) = \elr (A_p)$.
%\end{cor}

If a contractive $n$-tuple consists of isometries, then the invariants are both equal to the
wandering dimension of the $n$-tuple.

\begin{lem}\label{isom_invs}
If $S = (S_1, \ldots, S_n)$ is an $n$-tuple of isometries with pairwise orthogonal ranges, then
\[
\crv (S) = \elr (S) = \prk (S).
\]
\end{lem}

\Prf
Recall from Popescu's Wold decomposition that the $S_i$
are simultaneously unitarily equivalent to an orthogonal direct sum $T_i \oplus L_i^{(\alpha)}$,
where the $T_i$ form a representation of
the Cuntz algebra. The multiplicity $\alpha$ is equal to $\rk (I - \sum_{i=1}^n S_i S_i^*) = \prk
(S)$. The invariants are clearly stable under unitary equivalence. Thus, from the previous
lemma, $\crv (S) = \prk(S) \crv (L)$ and $\elr (S) = \prk (S) \elr (L)$.

However, a computation shows that
\[
\tr (Q_k) = \rk (Q_k) = 1 + n + \ldots + n^{k-1} = \frac{n^k - 1}{n-1}.
\]
Evidently, $\crv (L) = \elr (L) = \lim_{k\rightarrow\infty} \frac{n^k - 1}{n^k} = 1$, finishing
the proof.
\bx

The characterization of the pure rank from Lemma ~\ref{pure_rank},
together with the previous lemma, yields an esthetically pleasing result.

\begin{thm}\label{inv_ineq}
Let $A$ be a contractive $n$-tuple with minimal isometric dilation $S$. Then
\begin{itemize}
\item[$(i)$] $\crv (A) \leq \crv (S) = \prk (A)$
\item[and]
\item[$(ii)$] $\elr (A) \leq \elr (S) = \prk (A)$.
\end{itemize}
\end{thm}

\Prf
If $A=(A_1, \ldots, A_n)$ acts on $\H$ and $S=(S_1, \ldots, S_n)$ acts on a larger space $\K$,
then each $A_i$ is the compression to $\H$ of $S_i$. Further,  $\H$ sits inside $\K$ as a
co-invariant subspace for the isometries $S_i$. Thus,
\begin{eqnarray*} I_{\H} - \Phi_A^k(I_\H) &=& I_{\H} - \sum_{|w|=k} A_wA_w^* \\
                             &=& P_{\H} \Bigl( I_{\K} - \sum_{|w|=k}S_wS_w^* \Bigr) \mid_{\H}\\
               &=& P_{\H} ( I_{\K} - \Phi_S^k(I_\K)) \mid_{\H},
\end{eqnarray*}
which yields the inequalities. Further, from Lemma ~\ref{pure_rank},  the minimality of the
dilation guarantees that
\begin{eqnarray*}
\prk (S) &=& \rk (I_{\K} - \sum_{i=1}^n S_iS_i^* ) \\
&=& \rk (I_{\H} - \sum_{i=1}^n A_iA_i^*) =
\prk (A).
\end{eqnarray*}
An application of the previous lemma finishes the proof.
\bx

The following immediate corollary establishes the hierarchy of these invariants. It is exactly the
same as the analogous hierarchy from the commutative setting.

\begin{cor}\label{hierarchy}
Let $A$ be a contractive $n$-tuple of operators. Then
\[
0 \leq \crv (A) \leq \elr(A) \leq \prk(A).
\]
\end{cor}

The following sections contain an analysis of the basic properties of these functions, their
interplay and the information they can yield for a contraction.

%%%%%%%%%%%%%%%%%%%%
\section{Detection of Freeness}\label{S:basic_prop}
%%%%%%%%%%%%%%%%%%%%

In the commutative setting, the curvature invariant is sensitive enough to detect when a pure
contraction  is
`free'. That is, when a pure contraction is simultaneously unitarily equivalent to copies of the
multiplication operators by
the coordinate functions on symmetric Fock space. This is what Arveson calls its basic property.
In this section it is shown that the analogous result holds true for the non-commutative curvature
invariant. It can also be regarded as the basic property here. The notion of freeness for
non-commuting $n$-tuples can be expressed as proximity of an $n$-tuple to copies of the pure
isometries.

\begin{defn}
A contraction $A = (A_1, \ldots, A_n)$ is said to be {\bf free} if there is a positive integer
$\alpha$ for which $A$ is unitarily equivalent to $L^{(\alpha)} = (L_1^{(\alpha)}, \ldots,
L_n^{(\alpha)})$.
\end{defn}

The $\alpha$ above is of course necessarily the pure rank of $A$.
The proof of the detection of freeness theorem requires different methods in the
non-commutative setting.
The key technical device is the lemma which follows. It
depends on the existence of a limit related to the
curvature invariant. If $\M$ is a subspace of $\H_n$ which is invariant for $\fL_n$, then a
contractive $n$-tuple $A = (A_1, \ldots, A_n)$ is defined by
\[
A_i = P_{\M^\perp} L_i |_{\M^\perp} = (L_i^*|_{\M^\perp})^*.
\]
From the structure of the Frahzo-Bunce-Popescu dilation, all pure contractions are
direct sums of such $n$-tuples. By Theorem ~\ref{lim_exist} and Lemma ~\ref{isom_invs},
\[
1 = \crv (L) = \crv (A) + (n-1) \lim_{k\rightarrow\infty} \frac{\tr (P_\M (I - \Phi^k_L (I))
P_\M)}{n^k}.
\]
In particular, the latter limit exists.

\begin{defn}
Let $\M$ be a subspace of $\H_n$ which is invariant for $\fL_n$. Then define
\[
\ktil (\M) := (n - 1) \lim_{k\rightarrow\infty} \frac{\tr (P_{\M} (I - \Phi^k_L (I))
P_{\M})}{n^{k}}.
\]
\end{defn}

The advantage of considering this limit is that the lattice of $\fL_n$-invariant subspaces is well
known. These
subspaces are infinite dimensional if they are non-zero. The important observation here is that
the `wandering nature' of these subspaces forces the limit $\ktil (\M)$ to be non-zero exactly
when the subspace $\M$ is non-zero.

\begin{lem}\label{free_lem}
If $\M$ is a non-zero subspace of $\H_n$ which is invariant for $\fL_n$, then
$ \ktil (\M) > 0.$
\end{lem}

\Prf
From the decomposition theory for non-zero $\fL_n$-invariant subspaces, the subspace $\M$ has
a unique decomposition into cyclic invariant subspaces,
$ \M = \sum_{j} \oplus R_{\zeta_j} \H_n, $
where each $R_{\zeta_j}$ is an isometry in $\fR_n$ with $\fL_n$-wandering vector
$R_{\zeta_j} \xi_e = \zeta_j$. Thus,
\[
\tr (P_{\M}(I - \Phi^k_L (I))P_{\M}) = \tr(P_{\M} Q_k P_{\M} ) = \sum_{j} \tr(P_{j} Q_k
P_{j} ),
\]
where $P_j$ is the projection onto $R_{\zeta_j} \H_n$. Without loss of generality assume $\M =
R \H_n$, for an isometry $R = R_\zeta$ in $\fR_n$. Then the vectors $R\xi_w = R L_w \xi_e =
L_w \zeta$ form an orthonormal basis for $\M$.

Suppose $\zeta = \sum_{u\in\F_n} a_u \xi_u$ and let $v$ be a word of minimal length $|v| =
k_0$ such that $a_v = (\zeta , \xi_v) \neq 0.$ Then for $k> k_0$ and words $w$
with $|w| \geq k - k_0$, one has
\[
Q_k R\xi_w = Q_k L_w \zeta = Q_k ( I - Q_k) L_w \zeta  = 0.
\]
Whence, $(P_{\M}Q_k P_{\M} R\xi_w, R\xi_w ) = 0$. Conversely, for words $w$ with $|w| <
k - k_0$,
\begin{eqnarray*}
(Q_k R\xi_w, R\xi_w ) &=& \sum_{|u|\geq k_0} a_u (Q_k \xi_{wu}, L_w \zeta ) \\
                                      &=& \sum_{ |u| \geq k_0, \,\,\, |w| + |u| < k} a_u  (L_w \xi_u, L_w
\zeta) \\
                                      &=& \sum_{ |u| \geq k_0, \,\,\, |w| + |u| < k} |a_u|^2.
\end{eqnarray*}
In particular, for words $w$ with $|w| < k - k_0$, the lower bound
\[
(Q_k R\xi_w, R\xi_w ) \geq |a_v|^2
\]
is obtained. Hence, computing the trace yields
\[
\tr (P_{\M}Q_k P_{\M}) \geq |a_v|^2 (1 + n + \ldots + n^{k-k_0 -1}) = |a_v|^2
\frac{n^{k-k_0}-1}{n-1}.
\]
Therefore, it follows that $\ktil (\M) \geq \frac{|a_v|^2}{n^{k_0}} > 0.$
\bx

The non-commutative analogue can now be proved. As in the commutative setting, focusing on
the collection of pure contractions yields the best possible results here.

\begin{thm}\label{freeness}
If $A = (A_1, \ldots, A_n)$ is a non-zero pure contraction with $\prk(A) < \infty$, then
the following are equivalent:
\begin{itemize}
\item[$(i)$] $\crv(A) =  \prk (A)$
\item[$(ii)$] $A$ is free.
\end{itemize}
\end{thm}

\Prf
If $A$ is free, then the  invariants are equal by Lemma ~\ref{isom_invs}. Conversely, let
$\alpha = \prk(A)$ and suppose $A$ acts on $\H$. Note that since $A$ is
pure, its minimal isometric dilation is $L^{(\alpha)} = (L_1^{(\alpha)}, \ldots, L_n^{(\alpha)})$
acting on $\K := \H_n^{(\alpha)}$, which will be regarded as containing $\H$. Also,
 from the construction of the dilation the space $\H$  is co-invariant for the isometries
$L_i^{(\alpha)}$. Thus, since $P_{\H}(  I_{\K} - \Phi^k_{L^{(\alpha)}} (I))
P_{\H} = I_\H - \Phi^k_A (I)$ one has
\[
\tr (I_{\K} - \Phi^k_{L^{(\alpha)}} (I) ) =
                                                  \tr ( I_\H - \Phi^k_A (I)) + \tr (P_{\H^\perp}(I_{\K} -
\Phi^k_{L^{(\alpha)}} (I)) P_{\H^\perp}).
\]
Co-invariance also means that $\H^\perp$ is invariant for the algebra $\fL_n^{(\alpha)}$. Hence
$\H^\perp$ decomposes as a direct sum of $\alpha$ subspaces $\V_j$, each invariant for
$\fL_n$.
So the traces decompose as
\[
\tr (I_{\K} - \Phi^k_{L^{(\alpha)}} (I) ) =   \tr ( I_\H - \Phi^k_A (I)) + \sum_{j=1}^\alpha \tr
(P_{\V_j} Q_k P_{\V_j}).
\]

If $\H^\perp \neq \{0\}$, then some $\V_i \neq \{0\}$. Hence by Lemma ~\ref{free_lem},
one would have
%\begin{eqnarray*}
\[
\alpha = \crv (L^{(\alpha)}) \geq \crv(A) + \ktil (\V_i)
                                      > \crv (A) = \alpha ,
\]
%\end{eqnarray*}
an absurdity. This shows that $\H = \K$. Therefore, $A$ is its {\it own} minimal isometric
dilation
$L^{(\alpha)}$, and the result follows.
\bx

At the other end of the spectrum, there are large classes which are annihilated by the invariants.
This clearly includes contractions which are finite rank or of Cuntz type. Also, every
commuting contractive $n$-tuple is annihilated {\it exactly} because
Arveson's invariants exist.

\begin{prop}
If $A = (A_1, \ldots, A_n)$ is a contractive $n$-tuple of mutually commuting operators, then
$\crv(A) = \elr(A) = 0$.
\end{prop}

\Prf
The proof is the same for both invariants. For the curvature invariant,
\begin{eqnarray*}
\crv(A) &=& (n-1) \lim_{k\rightarrow\infty} \frac{\tr(I - \Phi^k(I))}{n^k} \\
              &=& (n-1)
\lim_{k\rightarrow\infty} \frac{k^n}{n^k} \frac{\tr(I - \Phi^k(I))}{k^n} = 0.
\end{eqnarray*}
\upbx

An interesting question which was posed by N. Spronk (University of Waterloo) is whether a
contraction is annihilated if a subset
of the $A_i$ are pairwise commuting.
For pure commuting contractions the proposition can also be observed by taking a spatial point
of view. For instance,
consider the motivating example there: $M = (M_1, \ldots, M_n)$ where $M_i$ is the
compression of $L_i$ to the symmetric Fock space $\H_n^s$. The basic point is that the words
of length at most $k$ in $n$ {\it commuting} letters is on the order of $\frac{k^n}{n!}$ for
large $k$.

Further consideration of this spatial point of view leads to a characterization of when the
curvature invariant annihilates a pure contraction. For the sake of brevity, the rest of this section
will focus on the $\prk(A) = 1$ case. A condition is obtained by Arveson in the commutative
pure setting in terms of sequences of {\it inner functions}. In the non-commutative setting,
typically isometries are taken to be {\it inner operators}. For instance, see
\cite{DP1,DP2,Pop_diln,Pop_fact,Pop_beur}
where a non-commutative inner-outer factorization theorem is obtained among other things.
Thus the following result can be thought of as an analogue of the commutative
characterization.

\begin{prop}
Let $A = (A_1, \ldots, A_n)$ be a pure contraction given by $A_i = P_{\M^\perp}
L_i|_{\M^\perp}$, where $\M$ is a subspace in $\Lat \fL_n$. Then the following are equivalent:
\begin{itemize}
\item[$(i)$] $\crv (A) = 0$,
\item[$(ii)$] There is a sequence of isometries $R_{\zeta_j}\in\fR_n$ with pairwise orthogonal
ranges for which $\{ \zeta_j \} \subseteq \M$ such that $\sum_j \ktil(R_{\zeta_j}\H_n) = 1$.
\end{itemize}
\end{prop}

\Prf
This is immediate from the definition of the limit $\ktil$ and the invariant subspace structure of
$\fL_n$. The first condition is satisfied exactly when $\ktil (\M) = 1$. However, when $\M \neq
\{ 0 \}$ it decomposes into an orthogonal direct sum $\M = \sum_j \oplus \M_j$ of cyclic
invariant subspaces $\M_j = R_{\zeta_j} \H_n = \overline{ \fL_n R_{\zeta_j}\xi_e}=
\overline{\fL_n \zeta_j}$. Hence, $\ktil (\M) = \sum_j \ktil(\M_j)$.
\bx

The rigidity of the rank function prevents a similar result for the Euler characteristic. However,
an interesting open problem is: Are the two invariants zero at the same time? This is the case in
all the classes which the author has come across. Roughly speaking, a pure contraction which is
annihilated by the curvature invariant can be thought of as lacking `freeness'. In other words,
such a contraction does not `look like' its dilation on any sustained patch. The following
example helps to illustrate this point further.

\begin{eg}
If $R \in \fR_n$ is a polynomial isometry of the form $R = \sum_{|w| =k}
a_w R_w$ for some $k \geq 0$ and $\sum_{|w| = k} |a_w|^2 = 1$, then a relatively simple
calculation shows that $\ktil (R\H_n) = \frac{n - 1}{n^k}$. Focusing on such isometries can
yield
interesting results. Indeed, they are part of the scenario of the two theorems which follow.  For
$k\geq 0$, the isometries $R_1^k R_2$ in $\fR_2$ have pairwise orthogonal ranges. They
determine a subspace in $\Lat \fL_2$ by $\M = \sum_{k \geq 0} \oplus \M_k$, where $\M_k =
R_1^k R_2 \H_2$. From the discussion above one has
\[
\ktil (\M) = \sum_{k\geq 0} \ktil (\M_k) = \sum_{k\geq 0} 2^{-k - 1} = 1.
\]
Thus the contraction $A = (A_1, A_2)$ defined by $A_i = P_{\M^\perp} L_i |_{\M^\perp}$
satisfies $\crv(A) = 0$. Notice that $A$ is infinite rank here; however, as the curvature invariant
suggests, it lacks freeness. In fact, $A$ is actually a {\it commuting} 2-tuple in this case. The
subspace $\M^\perp$ is given by
\[
\M^\perp = \spn \{ \xi_e, \xi_{1^k} : k\geq 1 \}.
\]
Hence it is clear that $A_1$ is unitarily equivalent to the unilateral shift, while $A_2 =
0$.
\end{eg}

A tight general characterization of when the invariants are in between the two extremes seems
unlikely. Instead, one must focus on particular classes of contractions. An excellent
example of how the curvature invariant can measure the `warping' of a contraction  is provided
in Example ~\ref{onedimldecay}.

A further investigation of the Euler characteristic allows one to view it as a measure of the
freeness of a pure contraction in certain cases. Recall that in the $\prk(A) = 1$ case  such a
contraction $A = (A_1,
\ldots, A_n)$, with each $A_i \in \BH$, has $L = (L_1, \ldots,
L_n)$ as its minimal isometric dilation and that $\H_n$ can be regarded as containing $\H$ as a
co-invariant subspace. This will be assumed to be the case for the next two theorems. To begin
with, a necessary condition on the size of both invariants is
obtained when $A$ possesses freeness.

\begin{thm}\label{elr1}
Let $A = (A_1, \ldots, A_n)$ be a pure contraction which acts on $\H$ with
$\prk(A) =1$. Suppose $\H^\perp$ is a cyclic $\fL_n$-invariant subspace of $\H_n$. If $\H$
contains $\spn \{ \xi_w :|w| \leq k \}$, then
\[
\elr (A) \geq \crv(A) > 1 - \frac{1}{n^k}.
\]
\end{thm}

\Prf
The inequality is trivial if $\H$ is all of $\H_n$. Otherwise, the subspace $\H^\perp$
can be written as $\H^\perp =  R_\zeta \H_n$ for an isometry $R_\zeta$ in $\fR_n$. The
generating wandering vector is $R_\zeta \xi_e = \zeta$, and an orthonormal basis for $\H^\perp$
is given by $ \{ L_w \zeta = R_\zeta \xi_w : w \in \F_n \}.$

Let $v$ be a word of minimal length $|v| = k_0$ such that $(\zeta , \xi_v) \neq 0$. The
hypothesis guarantees that $P_{\H^\perp} \xi_w =  P_{\H^\perp} Q_{k+1} \xi_w = 0$ for $|w|
\leq k$, so that $k_0 > k$.  For $l > k_0$ and words $w$ with $|w| \geq l - k_0$, one has
\[
Q_l R_\zeta \xi_w = Q_l L_w \zeta = Q_l (I - Q_{|w| + k_0}) L_w \zeta = 0.
\]
Hence, the range of $P_{\H^\perp}Q_l P_{\H^\perp}$ for $l > k_0$ is given by
\[
\ran (P_{\H^\perp}Q_l P_{\H^\perp}) = \spn \{ P_{\H^\perp}Q_l R \xi_w : |w| < l - k_0 \}.
\]
Thus, the following upper bound is obtained:
\begin{eqnarray*}
\tr (P_{\H^\perp}Q_l P_{\H^\perp}) &\leq& \rk (P_{\H^\perp}Q_l P_{\H^\perp}) \\
               &\leq& 1 + n + \ldots + n^{l-k_0 -1} = \frac{n^{l-k_0}-1}{n-1}.
\end{eqnarray*}
Therefore, $\ktil (\H^\perp) \leq 1 / n^{k_0}$. This finishes the proof since,
\[
1 = \crv(A) + \ktil (\H^\perp) \leq \crv(A) + \frac{1}{n^{k_0}} < \crv(A) + \frac{1}{n^k}.
\]
\upbx

This theorem does not hold when all subspaces $\H$ are considered. For example, when $\H$
is finite dimensional and contains $\spn \{ \xi_w :|w| \leq k \}$, both invariants are zero.
Further, a converse of this theorem does not hold for the curvature invariant, since the traces can
in general be spread over the entire space. However, the rigidity of the Euler characteristic can
be used to derive a related converse with more work.
In the cases for which the following theorem applies, the size of the Euler characteristic gives
direct information on what a contraction looks like in a sense.

\begin{thm}\label{elr2}
Let $A = (A_1, \ldots, A_n)$ be a pure contraction which acts on $\H$ with
$\prk(A) =1$. Suppose that $Q_l \H \subseteq \H$ for all sufficiently large $l$. Then the
subspace $\H$ of $\H_n$ contains $\spn \{ \xi_w :|w| \leq k \}$ when
\[
\elr(A) > 1 - \frac{1}{n^k}.
\]
\end{thm}

\Prf
Suppose $\H^\perp$ decomposes as $\H^\perp = \sum_{i=1}^\alpha \oplus R_{\zeta_i}\H_n$.
Let $k_0$ be a  minimal positive integer for which there is a word $v$ with $|v| = k_0$ and a $j$
such that
$(\zeta_j, \xi_v) \neq 0$. To prove the result, it suffices to show that $k_0 > k$. For if this is the
case, then $P_{\H^\perp} \xi_w = 0$ for words $w$ with $|w| \leq k$. Indeed, for such a word
$w$ and a typical basis vector $R_{\zeta_i} \xi_u$, one would have
%\begin{eqnarray*}
\[
(\xi_w, R_{\zeta_i} \xi_u) = (\xi_w, L_u \zeta_i)
               = (Q_{k+1} \xi_w, (I - Q_{k_0})L_u \zeta_i)
               = 0.
\]
%\end{eqnarray*}

The point here is that for all sufficiently large $l > k_0$, the vectors $Q_l R_{\zeta_j} \xi_w
=Q_l L_w \zeta_j$ for $0 \leq |w| < l- k_0$ are linearly independent. To see this, suppose $b_w$
are scalars for which the vector
$x = \sum_{0 \leq |w| < l - k_0} b_w Q_l L_w \zeta_j = 0$.
Since $l > k_0$, the vector $\xi_v$ satisfies $Q_l \xi_v = \xi_v$. Thus, by the minimality of the
word $v$,
\begin{eqnarray*}
0 = (x, \xi_v) &=& \sum_{0 \leq |w| < l - k_0} b_w (L_w \zeta_j, \xi_v) \\
          &=& b_e ( \zeta_j , \xi_v),
\end{eqnarray*}
so that $b_e = 0$. Now assume $b_w = 0$ for $0 \leq |w| < m < l-k_0$ and let $w_0$ be a word
of length $m$. Evidently, $Q_l \xi_{w_0 v} =  \xi_{w_0 v}$. Hence, again by minimality
\begin{eqnarray*}
0 = (x, \xi_{w_0 v}) &=& \sum_{0 \leq |w| < l - k_0}  b_w (L_w \zeta_j , \xi_{w_0 v}) \\
               &=& \sum_{|w| = m} b_w  (L_w \zeta_j , \xi_{w_0 v}) \\
                               &=& b_{w_0} (\zeta_j , \xi_v).
\end{eqnarray*}
Ergo, each $b_w = 0$. Thus for large $l$ such that $Q_l \H^\perp$ is contained in $\H^\perp$,
this shows that
$\rk (P_{\H^\perp}Q_l P_{\H^\perp}) \geq \frac{n^{l-k_0} - 1}{n-1}$.
Since the projection $Q_l$ reduces $\H$ for large $l$,
\[
\frac{n^l -1}{n-1} = \rk(Q_l) = \rk  (P_{\H}Q_l P_{\H}) + \rk (P_{\H^\perp}Q_l P_{\H^\perp}).
\]
It follows that
\begin{eqnarray*}
1 = \elr(A) + (n-1) \lim_{l\rightarrow\infty} \frac{\rk (P_{\H^\perp}Q_l P_{\H^\perp})}{n^l}
&\geq& \elr(A) + \frac{1}{n^{k_0}} \\
               &>& 1 - \frac{1}{n^k} + \frac{1}{n^{k_0}}.
\end{eqnarray*}
Therefore, $k_0 > k$ as required.
\bx

For pure contractions satisfying the conditions of this theorem, the conclusion quantifies the
statement: the closer $\elr (A)$ is to $\elr (L) = 1$, the more $A$ looks like $L$. Indeed, for
large $k$, the contraction $A$ is {\it equal} to $L$ on the large subspace $\spn \{ \xi_w : |w| < k
\}$.
Also notice that the limiting case as $k$ becomes arbitrarily large of the previous two theorems
is the statement of Theorem ~\ref{freeness}.

\begin{rem}
There are many examples which satisfy the hypothesis of the theorem.
Let $R$ be a polynomial isometry in $\fR_n$ of the form $R = \sum_{|v|=k} a_v R_v$ where
$\sum_{|v| = k} |a_v|^2 = 1$. Consider the cyclic
$\fL_n$-invariant subspace $R \H_n$. If $\H^\perp$ is the orthogonal direct sum of the ranges of
such isometries, then the contraction $A = (A_1, \ldots, A_n)$ acting on $\H$ by $A_i = P_{\H}
L_i |_\H$ fulfills the condition in the theorem. In fact, when the subspace $\H^\perp$ is cyclic
both of the previous theorems apply. As an example, consider $\H^\perp = R \H_n$ for the
isometry $R$ above. Then for $l > k$ and basis vectors $\xi_w$ one has
\[
Q_l R \xi_w = \sum_{|v| = k} a_v Q_l \xi_{wv} = \left\{ \begin{array}{cl}
          R \xi_w & \mbox{if $l > |w| + k$} \\
 0 & \mbox{otherwise.}
                         \end{array}
                    \right.
\]
In this case, the Euler characteristic can be computed directly as $\elr(A) = 1 - \frac{1}{n^k}$,
and $\H$ contains the subspace $\spn \{ \xi_w : |w| < k \}$.
\end{rem}

%The Euler characteristic can also be used to measure the freeness of more general decaying
%atomic representations. These representations are obtained by allowing for decaying outside of
%the central ring. For clarity, the $\prk(A) =1$ case will be looked at.
%
%\begin{eg}
%Consider the Hilbert space with orthonormal basis
%\[
%\{ \xi_{1,wi} : 2 \leq i \leq n \,\,\,\,{\rm and}\,\,\,\,w\in\F_n \}.
%\]
%Let $A = (A_1, \ldots, A_n)$ be a one-dimensional decaying $n$-tuple as defined in Example
%~\ref{decay}, with the proviso that the decaying occurs outside the ring. That is, there is a
%basis
%vector $\xi_{1,vi_0}$ with $|v| = l$ and a $j_0$ such that $A_{j_0} \xi_{1,vi_0}= \lambda
% \xi_{1,j_0vi_0}$ for some $|\lambda| < 1$. One still defines $A_j \xi_{1,vi_0}=
%\xi_{1,jvi_0}$
%for $j \neq i_0$. Also, $A_j \xi_{1,wi} = \xi_{1,jwi}$ for all other basis vectors and $1 \leq j
%\leq n$. The integer $l$ can be thought of as the {\it level of decaying}. In this case, it is easy
%to
%see that
%\[
%\ran (I - \Phi_A^k (I)) = \spn \{ \xi_{1,j_0vi_0} \}.
%\]
%
%Intuitively, such $n$-tuples seem to be `less free' when the decaying occurs further away from
%the ring. The Euler characteristic makes this idea concrete and provides a measure of the
%freeness of these representations. For the $n$-tuple $A$ above, a computation along the lines
%of
%Theorem ~\ref{decay_invs} shows that
%\[
%\elr (A) = \frac{n-1}{n^l}.
%\]
%Thus, $\elr (A)$ identifies where the decaying occurs.
%\end{eg}
%
%

%%%%%%%%%%%%%%%%%%%%%
\section{Continuity and Stability}\label{S:cont_stab}
%%%%%%%%%%%%%%%%%%%%%%%%%%%%

Any continuity result for these invariants must involve some sort of limit exchange, hence
uniform estimates are required. They are provided by the estimates obtained in Theorem
~\ref{lim_exist}. The abstract framework for proving semi-continuity comes from the
following elementary lemma.

\begin{lem}\label{tech_lemII}
Suppose $\{ s_k \}_{k\geq 1}$, $\{ a_k \}_{k\geq 1}$ and $\{ a^l_k \}_{k\geq 1}$ for $l \geq
1$ are sequences of non-negative real numbers which satisfy
\[ a_{k + 1} \leq a_k + s_k, \,\,\,\, a^l_{k + 1} \leq a^l_k + s_k \qfor l \geq 1 \qand \sum_{k \geq
1} s_k <
\infty.
\]
If $\lim_{k\rightarrow \infty} a_k = a$, $\lim_{k\rightarrow \infty} a^l_k = a^l$ for $l \geq 1$
and $\lim_{l\rightarrow \infty} a^l_k = a_k$ for $k \geq 1$, then
\[
\limsup_{l \geq 1} a^l \leq a.
\]
\end{lem}

\Prf
The proof is by contradiction. Without loss of generality, by dropping to a subsequence it can be
assumed that there is an $\eps > 0$ such that for $l \geq 1$,
\[
a^l - a \geq \eps > 0.
\]
Then choose an integer $K \geq 1$ for which $\sum_{m \geq K} s_m < \eps / 8$ and $|a_k - a| <
\eps / 4$ for $k \geq K$. Now fix $l_0 \geq 1$ and $k_0 \geq K$. It will be shown that the
sequence $\{ a_{k_0}^l \}_{l\geq 1}$ cannot converge to $a_{k_0}$.

There are two cases to consider. First, if $a_{k_0}^{l_0} \geq a^{l_0}$, then
\begin{eqnarray*}
| a_{k_0}^{l_0} - a_{k_0} | &\geq& | a_{k_0}^{l_0} - a| - | a_{k_0} - a| \\
               &>& (a_{k_0}^{l_0} - a^{l_0}) + (a^{l_0} - a) - \eps / 4 \geq 3\eps / 4.
\end{eqnarray*}
On the other hand, if $a_{k_0}^{l_0} < a^{l_0}$ one can choose $k > k_0$ for which
$a_{k_0}^{l_0} \leq a^{l_0}_k $ and $ | a_k^{l_0} - a^{l_0}| < \eps / 8$. From the
uniform estimates provided by the $s_k$,
\begin{eqnarray*}
0 \leq a_k^{l_0} - a_{k_0}^{l_0} &=& a_k^{l_0} - a_{k - 1}^{l_0} + \ldots + a_{k_0 +
1}^{l_0} - a_{k_0}^{l_0} \\
                    &\leq& \sum_{k_0 \leq m \leq k-1} s_k < \eps / 8.
\end{eqnarray*}
Hence,
\[
|a_{k_0}^{l_0} - a^{l_0}| \leq |a_{k_0}^{l_0} - a^{l_0}_k | + |a_k^{l_0} - a^{l_0}|< \eps / 4.
\]
Thus the following estimate is obtained:
\begin{eqnarray*}
|a_{k_0}^{l_0} - a_{k_0}| &\geq& |a_{k_0} - a^{l_0}| - |a_{k_0}^{l_0} - a^{l_0}| \\
               &>& |a - a^{l_0}| - |a_{k_0} - a | - \eps / 4 > \eps / 2.
\end{eqnarray*}
Therefore, $|a_{k_0}^{l_0} - a_{k_0}| > \eps / 2$ for $l_0 \geq 1$, and an absurdity is realized.
\bx

The notion of continuity here requires a spatial link with the associated contractions. Thus a
natural setting for considering the continuity of these invariants is with pure contractions,
because of the strong spatial link provided by the dilations. Recall that such a
contraction $A$ is the compression of $L^{(\alpha)}$ to an $(\fL_n^*)^{(\alpha)}$-invariant
subspace where $\alpha = \prk(A)$. For a pure contraction, let $P_A$ denote the projection onto
this determining co-invariant subspace.

\begin{thm}\label{proj_cont}
Suppose $A = (A_1, \ldots, A_n)$ and  $A_l = (A_{l,1}, \ldots, A_{l,n})$ for $l \geq 1$ are pure
contractions acting on the same space with
\[
\prk(A),  \prk(A_l)  \leq \alpha < \infty.
\]
If
$\wotlim_l P_{A_l} = P_A$, then
\[
\limsup_{l\geq 1} \crv(A_l) \leq \crv(A).
\]
\end{thm}

\Prf
The lemma is applied with $a_k^l = \tr (I - \Phi^k_{A_l}(I)) / n^k$, $a^l = \crv (A_l)$,
$a_k = \tr (I - \Phi^k_A (I)) / n^k$ and $a = \crv(A)$. The estimates are from Theorem
~\ref{lim_exist}:
\[
a_{k + 1} \leq a_k + \alpha / n^{k + 1} \qand a^l_{k + 1} \leq a^l_k + \alpha / n^{k + 1} \qfor l
\geq 1.
\]
Thus it suffices to show that $\lim_{l\rightarrow\infty} \tr (I - \Phi^k_{A_l}(I)) = \tr (I -
\Phi^k_A (I))$, for $k \geq 1$. This follows from weak convergence. For instance, when $\alpha
= 1$ the operators can be thought of as acting on $\H_n$, so that
\[
\tr (I - \Phi^k_{A_l}(I)) = \tr ( P_{A_l} Q_k P_{A_l}) = \tr ( Q_k P_{A_l}),
\]
for a fixed $k \geq 1$, and hence
\begin{eqnarray*}
\lim_{l \rightarrow \infty} \tr (I - \Phi^k_{A_l}(I)) &=& \lim_{l \rightarrow \infty} \sum_{|w| <
k} ( P_{A_l} \xi_w, \xi_w) \\
               &=& \sum_{|w| < k}  ( P_A \xi_w, \xi_w) = \tr ( (I - \Phi^k_A (I)).
\end{eqnarray*}
This completes the proof.
\bx

%\begin{note}
%In the infinite pure rank case there are a couple of possibilities. When $\tr (I - \Phi_A (I)) =
%\infty$ the theorem is trivial since $\crv (A) = \infty$. Otherwise, it is unclear at this point
%whether upper semi-continuity is satisfied.
%\end{note}

The upper semi-continuity of the Euler characteristic is addressed below. Simple examples
show that neither invariant is lower semi-continuous with respect to this convergence.

\begin{eg}
For $l \geq 1$ define contractions $A_l = (A_{l,1}, \ldots, A_{l,n})$ by $A_{l,i} = Q_l L_i
|_{Q_l \H_n} = (L^*_i |_{Q_l \H_n})^*$ for $1 \leq i \leq n$.  Then certainly
$\wotlim_{l\rightarrow \infty} Q_l = I$, so the associated
limit is $L = (L_1, \ldots, L_n)$. However, since $Q_lQ_kQ_l = Q_l$ for $k \geq l$ one has
\[
\rk ( I - \Phi^k_{A_l}(I)) = \rk ( Q_l Q_k Q_l ) = \frac{n^l - l}{n - 1}.
\]
Therefore,
\[
\crv (A_l) \leq \elr (A_l) = \lim_{k \rightarrow \infty} \frac{n^l - 1}{n^k} = 0.
\]
Whereas the limit satisfies $\crv(L) = \elr(L) = 1$.
\end{eg}

More general continuity results can be obtained for the curvature invariant. It turns out to be
upper semi-continuous with respect to coordinate-wise norm convergence. This can be proved
without focusing on the pure setting.

\begin{thm}\label{norm_semicont}
Suppose that $A = (A_1, \ldots, A_n)$ and $A_l = (A_{l,1}, \ldots, A_{l,n})$ for $l \geq 1$ are
contractions acting on the same space with
\[
\prk(A), \, \prk(A_l) \leq \alpha < \infty.
\]
If $\lim_{l\rightarrow\infty} || A_{l,i} - A_i || = 0$ for $1 \leq i \leq n$, then
\[
\limsup_l \crv(A_l) \leq \crv(A).
\]
\end{thm}

\Prf
As in Theorem ~\ref{proj_cont}, an application of Lemma ~\ref{tech_lemII} is the key here. In
particular, it suffices to check that $\lim_{l\rightarrow\infty} \tr (I - \Phi^k_{A_l}(I)) = \tr (I -
\Phi^k_A (I))$, for $k \geq 1$.

Fix $k \geq 1$. For ease of notation; let $B_l = I - \Phi^k_{A_l}(I)$ and  $B = I -
\Phi^k_A (I)$. Notice that
\begin{eqnarray*}
| \tr (B_l) - \tr (B) | &=& | \tr (B_l - B) | \\
               &\leq& || B_l - B || \rk (B_l - B).
\end{eqnarray*}
However, the ranks of $B$ and the $B_l$ are uniformly bounded above. Indeed, if $\Phi$ is one
of the
associated completely positive maps, then
\[
I - \Phi^k (I) = \sum_{i = 0}^{k - 1} \Phi^i ( I - \Phi(I)).
\]
It follows that, $\rk (B_l) \leq \alpha + n\alpha + \ldots + n^{k - 1}\alpha = \alpha \bigl(
\frac{n^k - 1}{n -
1} \bigr)$, for $l \geq 1$. The same is true for $B$. Hence one has,
\[
| \tr (B_l ) - \tr (B) | \leq \frac{2 \alpha (n^k - 1)}{n - 1}
 || \Phi^k_{A_l} (I) -  \Phi^k_A (I) ||,
\]
which converges to 0 as $l$ becomes arbitrarily large by hypothesis. Thus the result follows
from an application of the lemma.
\bx

This theorem is used in Section ~\ref{S:examples} to obtain interesting information
on the range of the curvature invariant. On the other hand,  the rigidity of the Euler characteristic
prevents non-trivial continuity results. The following example shows it is not upper
semi-continuous with respect to either of the natural notions of convergence.

\begin{eg}
Consider the sequences $\{ x_k \}_{k\geq 1}$ and $\{ y_k \}_{k\geq 1}$ of $\fL_3$-wandering
vectors belonging to $\H_3$ given by
\[
x_k = \alpha_k \xi_1 + \beta_k \xi_2 \qand y_k = \beta_k \xi_2 + \alpha_k \xi_3,
\]
where $\alpha_k, \beta_k > 0$, $\alpha_k^2 + \beta_k^2 = 1$ and $\lim_{k\rightarrow\infty}
\beta_k = 1$. Let $\M_k$ be the $\fL_3^*$-invariant subspace
\[
\M_k = \overline{\fL_3 x_k} \bigvee  \overline{\fL_3 y_k} \bigvee \spn \{ \xi_e \},
\]
with $P_k$ the projection onto $\M_k$. If $P$ is the projection onto $\spn \{ \xi_e, \xi_{u2} : u
\in \F_3 \}$, the first claim is that $\sotlim_{k\rightarrow\infty} P_k = P$.

To prove this, it suffices to check convergence on basis vectors $\xi_w$. This is trivial for
$\xi_e$. Consider a vector of
the form $\xi_{w1}$. Then $\xi_{w1}$ is orthogonal to all of the determining vectors of
$\M_k$ except $L_w x_k $ and $L_w y_k$. Let $z_k$ be the unit vector
orthogonal to $x_k$ obtained from the Gram-Schmidt process for which $\spn \{ x_k, y_k \} =
\spn \{ x_k, z_k \}$. It follows that $\{ L_w x_k, L_w z_k \}$ forms an orthonormal basis for
$\spn \{ L_w x_k, L_w y_k \}$. A computation shows that $z_k = \frac{1}{\sqrt{3}} ( -\alpha_k
\beta_k \xi_1 + \alpha_k^2 \xi_2 + \alpha_k \beta_k^{-1} \xi_3 )$. Hence,
\begin{eqnarray*}
P_k \xi_{w1} &=& ( \xi_{w1}, L_w x_k) L_w x_k + ( \xi_{w1}, L_w z_k) L_w z_k \\
          &=& \alpha_k L_w x_k  - \frac{\alpha_k \beta_k}{\sqrt{3}} L_w z_k,
\end{eqnarray*}
which converges to zero by hypothesis. An equivalent argument shows that
$\lim_{k\rightarrow\infty} P_k \xi_{w3} =  0$ for all words $w$ in $\F_3$. Next let
$\xi_{w2}$ be a basis vector in $P \H_3$. Then $\xi_{w2}$ is perpendicular to all the
determining vectors of $\M_k$ except $L_w x_k$ and $L_w y_k$. Another computation yields
\begin{eqnarray*}
P_k \xi_{w2} &=& ( \xi_{w2}, L_w x_k) L_w x_k + ( \xi_{w2}, L_w z_k) L_w z_k \\
          &=& \beta_k L_w x_k + \frac{\alpha_k^2}{\sqrt{3}} L_w z_k.
\end{eqnarray*}
Therefore, $\lim_{k\rightarrow\infty} P_k \xi_{w2} =  \xi_{w2}$. The first claim follows.

Hence the contractions $A = (A_1, \ldots, A_n)$ and $A_k = (A_{k,1}, \ldots, A_{k, n})$
defined
by
\[
A_i = P L_i P \qand A_{k,i} = P_k L_i P_k \qfor k\geq 1
\]
satisfy $\sotlim_{k\rightarrow\infty} P_k = P$. The second claim is that
\[
\lim_{k\rightarrow\infty} || A_{k,i} - A_i || = 0 \qfor 1 \leq i \leq n.
\]
To see this write $z_k = a_k\xi_1 + b_k\xi_2 + c_k\xi_3$
for $k\geq 1$. Then each of the coefficients belong to the  unit
disk and converge to zero as $k$ approaches infinity. If $x = \sum_{w\in\F_3}a_w \xi_w$ is a
unit vector in $\H_3$, then by symmetry and a  (long) computation, for $1 \leq i \leq n$ one
obtains
\begin{eqnarray*}
|| (A_{k,i} - A_i)x || &\leq& || a_e (A_{k,i} - A_i) \xi_e || + 2 || \sum_{w=u1} a_w (A_{k,i} -
A_i)\xi_w || \\
& & + ||  \sum_{w=u2} a_w (A_{k,i} - A_i)\xi_w || \\
     &\leq& \alpha_k \big( 1 + \frac{\beta_k^2}{3}\big)^{1 / 2} + 2\alpha_k
     (1 + \beta_k^2)^{1 / 2} \\
& & + \alpha_k \big[ \big( \beta_k + \frac{a_k \alpha_k}{\sqrt{3}} \big)^2
     + \big( \frac{\alpha_k b_k}{\sqrt{3}} - \alpha_k \big)^2
     + \frac{\alpha_k^2 c_k^2}{3}  \big]^{1  /  2}.
\end{eqnarray*}
In any event, this proves the claim since the upper bound converges to zero as $k$ becomes
arbitrarily large.

Thus this example satisfies the hypotheses of the previous two theorems. But it is easy to see that
\[
\rk (I - \Phi^l_{A_k}(I)) = 1 + 2 + 2\cdot 3 + \ldots + 2\cdot 3^{l -2} = 3^{l -1}
\]
and
\[
\rk (I - \Phi^l_A (I)) = 1 + 1 + 3 + \ldots + 3^{l - 2} = \frac{3^{l - 1} + 1}{2}.
\]
Therefore, $\elr (A_k) = 2 / 3$ for $k \geq 1$, while $\elr (A) = 1 / 3$. Hence the Euler
characteristic is not upper semi-continuous.
\end{eg}

This section finishes with a look at stability properties. It is obvious
that the invariants are stable under unitary equivalences. They are also stable under
multiplication by unitary matrices.

\begin{prop}
Let $A = (A_1, \ldots, A_n)$ be a contraction. If $U$ is an $n \times n$
unitary matrix, then $\crv (AU) = \crv (A)$ and $\elr (AU) = \elr (A)$.
\end{prop}

\Prf
Let $B = AU$. The point here is that the sequences of operators $\Phi_A^k (I)$ and $\Phi_B^k
(I)$
are the same. This is clear for $k=1$ since $\Phi_B(I) = BB^* = AA^* = \Phi_A(I)$.
For $k \geq 2$, let $\Phi_A^k (I)^{(n)}$ be the $n \times n$ diagonal matrix with $\Phi_A^k (I)$
down the
diagonal. Then by induction, one has
\begin{eqnarray*}
\Phi_B^k (I) &=& B ( \Phi_A^{k-1} (I)^{(n)} ) B^*
                    =  A U \Phi_A^{k-1} (I)^{(n)} U^* A^* \\
          &= & A \Phi_A^{k-1} (I)^{(n)} A^*
          =  \Phi_A^k (I).
\end{eqnarray*}
\upbx

It is natural to ask whether these invariants are stable under compressions to subspaces of finite
co-dimension. In certain situations they are.

\begin{prop}
Let $A = (A_1, \ldots, A_n)$ and $B = (B_1, \ldots, B_n)$ be contractions such that each $A_i$
is the compression of $B_i$ to a co-invariant subspace of finite
co-dimension. Then $\crv (A) = \crv (B)$ and $\elr (A) = \elr (B)$.
\end{prop}

\Prf
Suppose $B$ acts on $\H$ and $A$ acts on a subspace $H_0$ for which $\H_1 := \H \ominus
\H_0$ is finite dimensional. By hypothesis, each $A_i = P_{\H_0} B_i |_{\H_0} =
(B_i^* |_{\H_0})^*$, whence $\Phi^k_A (I) = P_{\H_0} \Phi^k_B (I) |_{\H_0}$ for $k \geq 1$.
Let $B_k = I - \Phi^k_B (I)$. It follows that
\begin{eqnarray*}
\rk ( I - \Phi^k_A (I)) &\leq& \rk (B_k) \\
          &\leq& \rk (I - \Phi^k_A (I)) +\rk (P_{\H_0} B_k P_{\H_1}) + \rk ( P_{\H_1}
B_k).
\end{eqnarray*}
Thus, the associated Euler characteristics are evaluated as
\begin{eqnarray*}
|\elr(B) - \elr(A)| &=& (n-1) \lim_{k\rightarrow\infty} \Big( \frac{\rk (B_k)}{n^k} -
\frac{\rk (I - \Phi^k_A (I))}{n^k} \Big) \\
                             &\leq& (n-1) \lim_{k\rightarrow\infty}\frac{2 \dim
\H_1}{n^k} = 0.
\end{eqnarray*}
An easier computation works for the curvature invariant since the trace is linear.
\bx

Arveson's invariants are stable under compressions to {\it any} subspace of finite co-dimension.
This
is not true for the non-commutative versions, even if the subspace is invariant for the $n$-tuple.

\begin{eg}
For $1 \leq i \leq n$, let $A_i$ be the compression of $L_i$ to the $\fL_n$-invariant subspace
$\xi_e^\perp := \spn \{ \xi_w : |w| \geq 1 \}$. Then $A$ and $L$ fit into the context of the
previous proposition. Recall that $\crv(L) = \elr(L) = 1$. On the other hand a simple calculation
shows that $\prk(A) = \rk ( I - \Phi_A(I)) = n$. In fact, the
contractive $n$-tuple $A = (A_1, \ldots, A_n)$ satisfies
\begin{eqnarray*}
\Phi_A^k (I) \xi_w =  \sum_{|v| =k} A_v A_v^* \xi_w &=& \sum_{|v| =k} L_v P_{\xi_e^\perp}
L_v^* \xi_w \\
          &=& \left\{ \begin{array}{cl}
                                                            0 & \mbox{if $1 \leq |w| \leq k$} \\
                         \xi_w & \mbox{if $|w| > k$.}
                         \end{array}
                    \right.
\end{eqnarray*}
Thus by computation one has,
\[
\crv(A) = \elr(A) = (n-1) \lim_{k\rightarrow\infty} \left[ \frac{n+n^2+ \ldots +
n^k}{n^k}\right] = n.
\]
It follows from Theorem ~\ref{freeness} that $A \simeq L^{(n)}$. This can also be
observed by
noting that $\xi_e^\perp$ decomposes as the orthogonal direct sum of $n$ subspaces which
reduce $A$.
The compression of $A$ to each of these subspaces is unitarily equivalent to $L$. Of course,
this example doesn't work in the commutative setting. The reason is that the $A_i$ would act
on symmetric Fock space where the version of $\xi_e^\perp$ does not decompose into a direct
sum in this manner because the $n$ subspaces have overlap.
\end{eg}

%%%%%%%%%%%%%%%%%
\section{Examples}\label{S:examples}
%%%%%%%%%%%%%%%%%

It is probably not reasonable to expect a tight characterization of pure contractive $n$-tuples,
since they are the multi-variable analogues of completely non-unitary operators. Nonetheless,
this
section contains a rich collection of examples for which $I - \Phi(I)$ is finite rank. In particular,
a new class of examples is introduced which fills out the range of the curvature invariant and
illustrates further how the curvature of an $n$-tuple can be measured. This class consists of
finite rank perturbations of certain Cuntz representations which have an overlap with wavelet
theory. Their
Euler characteristic and curvature invariant can be computed directly, hence information on the
general relationship between the two is obtained. There is a whole host of open problems
regarding the ranges of the invariants, and some of these are pointed out below.

\begin{eg}\label{decay}
In the paper \cite{DP1},  Davidson and  Pitts described a class of representations of the
Cuntz-Toeplitz
$\ca$-algebra $\E_n$ which they called atomic free semigroup representations. These
representations
decompose as a direct integral of irreducible atomic representations, which are of three types.
The first is the left regular representation, and is the only one which does not factor through the
Cuntz algebra $\O_n$. The second type is a class of inductive limits of the left regular
representation which are classified by an infinite word up to tail equivalence. The third type
can be called atomic ring representations. These representations have the shape of a benzene
ring, with infinite trees leaving each node. The nodes correspond to basis vectors and each tree
swept out is a copy of the left regular representation. The associated isometries with pairwise
orthogonal ranges map ring basis vectors either to the next vector in the ring, allowing for
modulus one multiples of the image vector, or to a top of the tree which lies below the
original node.

These representations can be perturbed to obtain new examples which fit into the context of
this paper. The idea is to preserve the structure of the ring representations, with the proviso that
the images of vectors lying in the ring are allowed to be strictly contractive multiples,
instead of just modulus one multiples. These new representations can be thought of as
possessing a certain decaying property as one moves around the ring.

The construction proceeds as follows: Suppose $u = i_1i_2\ldots i_d$ is a word in $\F_n$. Let
$\H_u$ be the Hilbert space with orthonormal basis,
\[
\{ \xi_{s,w} : 1 \leq s \leq d {\rm \,\,\,and\,\,\,} w \in \F_n \setminus \F_n i_s \}.
\]
If $\vec{\lambda}= (\lambda_1, \ldots , \lambda_d)$ is a $d$-tuple of complex scalars for which
$|| \vec{\lambda} ||_{\infty}\leq 1$, then define a contraction $A = (A_1, \ldots , A_n)$
acting on
$\H_u$ by
\[
\begin{array}{rclcl}
A_i \xi_{s,e} &=& \lambda_s \xi_{s + 1,e} & {\rm if} & i=i_s, \,1\leq s \leq d \\
A_i \xi_{s,e} &=& \xi_{s,i} & {\rm if} & i \neq i_s \\
A_i \xi_{s,w} &=& \xi_{s,iw} & {\rm if} & w \neq e
\end{array}
\]
So the ring vectors are given by $\xi_{s,e}$ for $1\leq s \leq d$. The associated representation of
$\F_n$ is denoted by $\sigma_{u,\vec{\lambda}}$. When each $\lambda_s$ is modulus one, the
resulting representation is a Cuntz representation in the class described by Davidson and Pitts.
However, when some $\lambda_s$ is on the {\it open} unit
disk, there is actual decaying which occurs around the ring. Hence the associated representation
will be called a {\bf decaying atomic representation}. The dimension of the subspace determined
by the central ring vectors is referred to as the dimension of the representation.
\end{eg}

The one-dimensional decaying contractions are in fact finite rank perturbations of the Cuntz
representation which gives rise to the Haar basis wavelet. This example is analyzed further
below. In
general the pure rank of these representations is determined by the amount of decaying
which occurs.

\begin{prop}
If $\sigma_{u,\vec{\lambda}}$ is the decaying atomic representation associated with the word
$u = i_1 \cdots i_d$ and vector $\vec{\lambda} = (\lambda_1, \ldots, \lambda_d)$, then the rank
of $I - \Phi (I)$ is equal to the cardinality of the set
$ \{ s : | \lambda_s | < 1 \}.$
\end{prop}

\Prf
This is straight from a computation of $I - \Phi(I)$ on the determining basis for the
representation. For ring basis vectors $\xi_{s, e}$ one has $A_i^* \xi_{s, e} = 0$ when $i \neq
i_{s-1}$. Whence
\[
(I - \Phi(I) ) \xi_{s, e} = (I - A_{i_{s-1}}A^*_{i_{s-1}}) \xi_{s, e} = (1 - | \lambda_{s-1} |^2)
\xi_{s, e}.
\]

On the other hand, every basis vector outside the ring is of the form $\xi_{s, wi}$ for some
 letter $i_s$ in $u$ and $i \neq i_s$. If $w i = j_0 v$, where $1 \leq
j_0 \leq n$ and $v$ is a word in $\F_n$, then $A_j^* \xi_{s, wi} = 0$ for $j \neq j_0$. Thus,
\[
(I - \Phi(I) ) \xi_{s, w i} = (I -  A_{j_0} A^*_{j_0}) \xi_{s,j_0v} = 0.
\]
\upbx

These examples form a tractable class of pure contractions. This also shows how a
Cuntz $n$-tuple can be perturbed by a finite rank operator to become pure.

\begin{thm}
Let $A = (A_1, \ldots, A_n)$ be a decaying atomic contraction. Then A is a pure
contraction.
\end{thm}

\Prf
Suppose $A$ is determined by a word $u = i_1 \cdots i_d$ and vector $\vec{\lambda} =
(\lambda_1, \ldots, \lambda_d)$. It is required to show that $\Phi^\infty (I) = 0$. Equivalently,
$\lim_{k\rightarrow\infty} \Phi^k(I) \xi_{s,w} = 0$, for all basis vectors $\xi_{s,w}$.

Consider a fixed basis vector $\xi_{s,w}$ where $1 \leq s \leq d$ and $w \in \F_n \setminus \F_n
i_s$. For a given $k$, there is only one word $v_k$ of length $k$ for which $A_{v_k}^*
\xi_{s,w} \neq 0$. For sufficiently large $k$, this word $v_k$ will pull $\xi_{s,w}$ back toward
the benzene ring, and then move around the ring. Such a word will be of the form
\[
v_k = w i_{s-1}\cdots i_d u^{l_k} i_1 \cdots i_{m_k}, \qfor \quad\text{some}\quad
1 \leq m \leq d.
\]
Clearly $l_k$ becomes arbitrarily large as $k$ does. Let $r = \max \{|\lambda_j| : |\lambda_j|
< 1\}$.
Then,
\begin{eqnarray*}
|| \Phi^k(I) \xi_{s,w} || & = & || \sum_{|v| = k} A_v A_v^* \xi_{s,w}||
                                       = || A_{v_k} A_{v_k}^* \xi_{s,w}|| \\
                                       &=& || A_{v_k} A_{i_{m_k}}^* \cdots A_{i_1}^* (A_u^*)^{l_k}
\xi_{d,e} ||    \leq || (A_u ^*) ^{l_k} \xi_{d, e} ||   \leq r^{l_k}.
\end{eqnarray*}
Hence, $\lim_{k\rightarrow\infty} || \Phi^k (I) \xi_{s, w}|| = 0$ as required.
\bx

In principle, the Euler characteristic and curvature invariant of these decaying atomic
contractions can be computed directly. The general pure rank one proof of the Euler
characteristic is included first.

\begin{lem}\label{decay_invs}
Let $A = (A_1, \ldots ,A_n)$ be a $d$-dimensional decaying atomic contraction with $\prk(A)
=1$
determined by a scalar $\lambda$ with $0\leq | \lambda | < 1$. Then
\[
\elr(A) = 1 - \frac{1}{n^d}.
\]
\end{lem}

\Prf
Without loss of generality assume the $n$-tuple $A$ is determined by the representation
$\sigma_{u, \vec{\lambda}}$ where $u = i_1\cdots i_d$ and $\vec{\lambda}$ is
the $d$-tuple $\vec{\lambda} = (\lambda, 1, \ldots ,1)$. The associated orthonormal basis is
\[
\{ \xi_{s,w} : 1 \leq s \leq d {\rm \,\,\,and\,\,\,} w \in \F_n \setminus \F_n i_s \}.
\]
Let $r = |\lambda |^2$ and let $R_k = I - \Phi^k(I) = I - \sum_{|v|=k}A_vA_v^*$. The action
around the ring is given by $A_{i_1} \xi_{1,e} = \lambda \xi_{2,e}$ and $A_{i_s} \xi_{s,e} =
 \xi_{s + 1,e}$ for $2 \leq s \leq d$ (where $d + 1$ is identified with $1$).

Consider a typical basis vector $\xi_{s,e}$ with $2 \leq s \leq d+1$ and put $m = s - 2 \geq 0$.
Let
$k \geq d$ be a positive integer and let $w$ be a word in $\F_n$ with $|w| = k - s +1$. Then
\begin{eqnarray*}
R_k \xi_{s,wi} &=& \xi_{s,wi} - A_{wi i_{s-1}\cdots i_2} A_{i_2}^* \cdots A^*_{i_{s-1}}
A_i^* A_w^*  \xi_{s,wi} \\
                        &=&  \xi_{s,wi} - A_{wi i_{s-1}\cdots i_2} \xi_{2,e} = 0.
\end{eqnarray*}
Similarly, $R_k \xi_{s,wi} = 0$ for $|w| \geq k - s +1$.

On the other hand let $w$ be a word with $|w| = k - s$. Then
\begin{eqnarray*}
R_k \xi_{s,wi} &=& \xi_{s,wi} - A_{wi i_{s-1}\cdots i_1} ( A_{wi i_{s-1}\cdots i_1})^*
\xi_{s,wi} \\
                        &=& ( 1 - r ) \xi_{s,wi}.
\end{eqnarray*}
Analogously, for $k \geq d$, every basis vector $\xi_{s,wi}$ with $|w| \leq k - s$ will belong to
the range of $R_k$ since $R_k \xi_{s,wi} = (1 - r^t) \xi_{s,wi}$, for some $t$ depending on
$|w|$ and $d$. The total number of such vectors is
\[ 1 + (n-1) + (n-1)n + \ldots + (n-1)n^{k-s} = n^{k - s + 1}. \]
Therefore the rank of $R_k$ for $k \geq d$ is computed as
\[
\rk (R_k) = n^{k-1} +  n^{k-2} + \ldots + n^{k-d} = \frac{n^k}{n-1} \bigl( 1 - \frac{1}{n^d}
\bigr),
\]
which shows that $\elr(A) = 1 - \frac{1}{n^d}$, as claimed.
\bx

The formula for the curvature invariant of the one-dimensional $n$-tuples is readily obtained
from the analysis in the previous proof.

\begin{lem}
Let $A = (A_1, \ldots, A_n)$ be a one-dimensional decaying atomic contraction determined by a
scalar $\lambda$ with $0 \leq |\lambda| <1$ (hence $\prk(A) =1$). Then
\[
\crv (A) = (n-1) \frac{1-| \lambda |^2}{n- | \lambda |^2}.
\]
\end{lem}

\Prf
Without loss of generality assume $u=1$ so that the central ring vector is $\xi_{1,e}$, and
define $r$ and $R_k$ as above. Then the ring action is given by  $A_1^* \xi_{1,e} =
\bar{\lambda} \xi_{1,e}$ and $A_i^* \xi_{1,e}= 0$ for $i \neq 1$. Let $k \geq 2$ be a fixed
positive integer. Then
%\begin{eqnarray*}
\[
R_k \xi_{1,e} =  \xi_{1,e} - A_1^k (A_1^*)^k  \xi_{1,e}
                        =   (1- r^k) \xi_{1,e}.
\]
%\end{eqnarray*}
Further, if $|w| = k - l$ for some $2 \leq l \leq k$, then
%\begin{eqnarray*}
\[
R_k \xi_{1, wi} = \xi_{1, wi} - A_w A_i A_1^{l-1} (A_1^*)^{l-1} A_i^* A_w^*
\xi_{1,wi}
                           = (1 - r^{l-1}) \xi_{1, wi}.
\]
%\end{eqnarray*}
However, if $w$ is a word of length at least $k-1$, say $wi = uv$ with $|u| =k$, then
\[
R_k \xi_{1, wi} =  \xi_{1, wi} - A_u A_u^* \xi_{1, uv} = 0.
\]

Therefore the traces can  be evaluated as
\begin{eqnarray*}
\tr (R_k) &=& 1\! - \! r^k + (n \! - \!1) \big[ (1\! - \! r^{k-1}) + n (1\! - \! r^{k-2}) + \ldots +
n^{k-2}(1\! - \! r) \big] \\
               &=& n^{k-1} - r^k -(n-1) n^{k-1} \left[ \frac{r}{n} \frac{(r/n)^k - 1}{r/n - 1}\right] \\
               &=& n^{k-1} - r^k - \frac{(n-1)r}{n(n-r)} (n^k - r^k).
\end{eqnarray*}
Thus the curvature invariant is given by
\begin{eqnarray*}
\crv (A) &=& (n-1) \lim_{k\rightarrow\infty} \left[ \frac{1}{n} - \frac{r^k}{n^k} -
\frac{(n-1)r}{n(n-r)} (1 - \frac{r^k}{n^k}) \right] \\
              &=& \frac{n-1}{n} \left[ 1- \frac{(n-1)r}{n-r} \right] \\
              &=& (n-1) \frac{1-r}{n-r},
\end{eqnarray*}
as claimed.
\bx

\begin{rem}
In the commutative setting the invariants are always integers. In fact
for the  examples focused on by Arveson in \cite{Arv_euler} (those for which the
associated module is graded), the curvature invariant
and Euler characteristic are always equal.
This is not the case for the non-commutative versions. Indeed, the decaying atomics provide an
interesting tractable class of examples for which the invariants are distinct.
For instance, even consider the one dimensional decaying atomic $2$-tuple associated with
$\lambda = 1/ \sqrt{2}$. The theorem tells us that $\crv (A) = 1 / 3$ and $\elr (A) = 1 / 2$ in this
case. This class is analyzed further in the example which follows below. At this point there does
not appear to be a good general characterization of when
the two invariants are equal.  It is possible to say things in special cases. Indeed,
it appears that for a decaying contraction $A = (A_1, \ldots, A_n)$ determined by the
representation
$\sigma_{u,\vec{\lambda}}$, the condition $\crv(A) = \elr(A)$ is satisfied exactly when one of
the two extreme cases occurs. That is, the vector $\vec{\lambda}$ is either
$\vec{\lambda} = (0, \ldots, 0)$ or $\vec{\lambda} = (1,\ldots, 1);$ in other words, either there
is full annihilation around the ring, or there is no decaying at all and therefore
$\sigma_{u,\vec{\lambda}}$ is a Cuntz representation. Along these lines, there may be a
non-commutative Gauss-Bonnet-Chern theorem which is the analogue of the commutative
operator-theoretic version from \cite{Arv_euler}.
\end{rem}

The decaying atomics are all pure contractions; however, it is difficult to see in general exactly
how they sit inside their dilation. For the one-dimensional contractions this is both possible and
beneficial as it relates to the curvature invariant.

\begin{eg}\label{onedimldecay}
A one-dimensional decaying atomic  $n$-tuple is characterized by (without loss of generality)
$A_1 \xi_{1,e} = \lambda \xi_{1,e}$, for some $|\lambda| < 1$ and $A_i \xi_{1,e} = \xi_{1,i}$
for $i \neq 1$. The contractions act as the left regular representation below the ring, or  `loop' in
this case. Now, each
$A_i = P_\M L_i |_\M = (L_i^*|_\M)^*$ is the compression of $L_i$ to a co-invariant subspace
$\M$ of $\H_n$. Notice that one has
\[
L_1^* \xi_{1,e} = A_1^* \xi_{1,e} = \bar{\lambda} \xi_{1,e} \qand L_i^* \xi_{1,e} = A_i^*
\xi_{1,e} = 0 \qfor i \neq 1.
\]
Thus, $\xi_{1,e}$ is an {\it eigenvector} of the algebra $\fL_n^*$. But these vectors were
completely described in \cite{DP1}. In particular, $\xi_{1,e}$ is the eigenvector corresponding
to the scalar $n$-tuple $\vec{\lambda} = (\lambda,0,\ldots, 0)$. That is,
\[
\xi_{1,e} = \nu_\lambda := \sqrt{1 - |\lambda|^2} \sum_{k\geq 0} \bar{\lambda}^k \xi_{1^k}.
\]
It follows that the subspace $\M$ is given by
\[
\M = \spn \{ \nu_\lambda, L_w \nu_\lambda : w \in\F_n\setminus \F_n1 \}.
\]

From the previous lemma, the curvature invariant of the contraction is given by
\[
\crv(A) = (n-1) \frac{1 - |\lambda|^2}{n - |\lambda|^2}.
\]
For $\lambda = 0$, there is no loop and $\nu_0 = \xi_e$ is the vacuum vector. Thus $A$ sits
{\it cleanly} inside its dilation relative to the spatial structure of the underlying Fock space.
However, as
$|\lambda|$ becomes positive and increases toward 1, one sees that $\nu_\lambda$ spreads over
the whole space and $A$ becomes more `warped' or `curved' as it sits inside the dilation. The
curvature invariant reflects this line of thinking. Indeed, it decreases as $|\lambda|$ increases and
approaches zero as the associated representation warps right into a Cuntz representation.
\end{eg}

The rest of this section focuses on the ranges of the invariants. Some of the open questions are
pointed out. The variety of even the one-dimensional decaying $2$-tuples turns out to be
extensive enough to obtain the entire positive real line in the image of the curvature invariant.

\begin{thm}
For every $r \geq 0$, there is a contraction $A = (A_1, A_2)$ for which $\crv (A) = r$.
\end{thm}

\Prf
By Lemma ~\ref{dir_sum}, it is sufficient to obtain an interval in the range of $\crv (A)$ which
includes 0. This is just a matter of using the previous lemma and solving an identity. For positive
numbers $r$ with $0 \leq r \leq 1/ 2$, the number $s = \frac{1-2r}{1-r}$ belongs to the unit
interval. The one dimensional decaying $2$-tuple $A$ determined by $\lambda = \sqrt{s}$
satisfies $\crv (A) = r$.
\bx

\begin{rem}
There is not as much information available on the range of the Euler characteristic.
It is easy to construct examples $A$  which satisfy $\elr(A) = 1 / n$. Indeed, defining $A$ by
$A_i = P_\M L_i |_\M$ where $\M^\perp = \sum_{i=2}^n \oplus R_i \H_n$ suffices. Thus, using
direct sums, it follows that every positive rational number is in the
range of the Euler characteristic. It seems reasonable to make the guess that every positive real is
in the range of $\elr(A)$,  it would be surprising if this were not the case.
\end{rem}

\begin{note}\label{cuteeg}
As mentioned in the introduction, G. Popescu has shown that the range of the Euler
characteristic is indeed the positive real line. There is another cute way to see the full ranges of
the invariants, which was obtained by the author after submission. For $k \geq 0$, define an
$\fL_2^*$-invariant subspace by
\[
\M_k = R_1^k R_2 \H_2 \oplus \spn \{ \xi_{1^j} : 0 \leq j \leq k \}.
\]
Every $r$ in the unit interval can be represented by a binary expansion $r = \sum_{k\geq 0}
\eps_k
2^{-k-1}$ with each $\eps_k \in \{ 0,1\}$. Let $\eps_k \M_k$ be $\{0\}$ if $\eps_k = 0$ or
$\M_k$ if $\eps_k = 1$. Let $\M_r$ be the closed span of the subspaces $\eps_k \M_k$ for $k
\geq 0$. Then a row contraction $A_r = (A_{1,r}, A_{2,r})$ is defined by
$A_{i,r} = P_{\M_r} L_i|_{\M_r} = (L_i^*|_{\M_r})^*$. A computation shows that
$\crv (A_r) = \elr (A_r) = r$. Thus verifying that the unit interval (and hence the positive real
line) belongs to the ranges of both invariants. This example is quite satisfying as it is easy to
grasp onto and since it realizes the full ranges of the invariants with row contractions for
which the two are equal.
\end{note}

For decaying contractions with higher dimensional central rings the formulae for the curvature
invariant become particularly nasty. Nonetheless, the continuity results from the previous section
can allow one to avoid these computations and obtain fruitful results.

\begin{thm}
For every $\eps > 0$ and positive integers $k \geq 1$ and $ n \geq 2$, there is a contraction $A =
(A_1, \ldots, A_n)$
for which $\prk(A) = k$,
\[
0 < \crv(A) < \eps \qand  \elr(A) > k - \eps.
\]
\end{thm}

\Prf
It is sufficient to prove the $\prk(A) =1$ case since direct sums can then be used in the general
case. The pure rank one $d$-dimensional decaying atomic $n$-tuples provide the concrete
examples here.
The Euler characteristic is always $\elr(A) = 1 - \frac{1}{n^d}$, independent of the decaying
factor $\lambda$. Hence by choosing large enough central rings, $\elr(A)$ asymptotically
approaches $1$.

Given a fixed word $u$ in $\F_n$ of length $d$, consider the decaying atomic $n$-tuple $A_r =
(A_{1,r}, \ldots, A_{n,r})$ acting on $\H_u$ which is determined by the $d$-tuple $\vec{r} =
(r,1, \ldots, 1)$ for $0 \leq r \leq 1$. When $r=1$ the $n$-tuple forms a Cuntz representation, so
that $\crv(A_1) =0$. However, observe that if $x = \sum_{s,w,i}a_{s,w,i} \xi_{s,wi}$ is a unit
vector in $\H_u$, then
\begin{eqnarray*}
|| (A_{j,1} - A_{j,r})x|| &=& || \sum_{s,w,i} a_{s,w,i} (A_{j,1} - A_{j,r}) \xi_{s,wi} || \\
               &=& || a_{1,e,e} (A_{j,1} - A_{j,r}) \xi_{1,e} || \\
               &=& \left\{ \begin{array}{cl}
          | a_{1,e,e} | | 1 - r | & \mbox{if $j = i_1$} \\
 0 & \mbox{if $j \neq i_1$}
                         \end{array}
                    \right.
\end{eqnarray*}
It follows that $\lim_{r\uparrow 1} || A_{j,1} - A_{j,r} || = 0$ for $1 \leq j \leq n$. Hence
$\limsup_r \crv(A_r) = \crv(A_1) = 0$  by the
upper semi-continuity of $\crv (\cdot)$ proved in Theorem ~\ref{norm_semicont}.
This finishes the proof.
\bx

\begin{rem}
Thus the invariants can be asymptotically as far apart as possible. It would be interesting to
know if the extreme case can be attained. In other words, is there an $A$ with $\prk(A) =1$
such that $\crv(A) = 0$ and $\elr(A) = 1$?  This relates to the earlier question of whether the two
invariants always annihilate at the same time.
\end{rem}

There is another class of examples which in a sense are pervasive. If $\M$ is a subspace of
$\H_n$ which is co-invariant for $L = (L_1, \ldots, L_n)$, then a contraction $A =
(A_1, \ldots, A_n)$ is defined by $A_i = P_\M L_i |_\M = (L_i^*|_\M)^*$ for $1 \leq i \leq
n$. Recall that co-invariance shows
%\begin{eqnarray*}
\[
I - \Phi_A^k(I) = I_\M - \sum_{|w| = k}A_w A_w^*
          = P_\M ( I - \sum_{|w| = k}L_w L_w^*) |_\M.
\]
%\end{eqnarray*}
Thus, $L$ is an isometric dilation of $A$ which is minimal when $\xi_e$ belongs to $\M$.
Also recall that from the structure of the Frahzo-Bunce-Popescu dilation,  all pure
contractions can be obtained from direct sums of such $n$-tuples. It turns out that
examples can be constructed from this point of view which fill out the range of the curvature
invariant. The proof uses the invariant defined in Section ~\ref{S:basic_prop}.

\begin{thm}
For every $r \geq 0$, there are positive integers $n \geq 2$ and $k \geq 1$ and a subspace $\M$
in $\Lat (
\fL_n^*)^{(k)}$ for which the contraction $A = (A_1, \ldots, A_n)$ defined by $A_i = P_{\M}
L_i^{(k)} |_{\M}$ satisfies $\crv(A) = r$.
\end{thm}

\Prf
The case $\crv(A) =0$ is covered by $\M = \{ 0 \}$. Recall from the remarks preceding
Lemma ~\ref{free_lem} that if $\M$ is an $\fL_n^*$-invariant subspace, then
$
1 = \crv(L) = \crv(A) + \ktil (\M^\perp).
$
Thus it suffices to capture the unit interval in the range of $\ktil$.

Consider $0 < r \leq 1 / 4$. Choose $n \geq 3$ such that $1/ n^2 < r \leq 1 / (n-1)^2$. A
computation shows that $1 - nr > 0$ and $\frac{n}{n-1} (1-nr) < 1$. Let
\[
a_2 = \sqrt{\frac{n}{n-1}(1-nr)} \,\,\,\, {\rm and} \,\,\,\, a_1 = \sqrt{1 - a_2^2}.
\]
Define an isometry $R$ in $\fR_n$ by $R = a_1 R_1 + a_2 R_2^2$. Let $\M^\perp = R\H$.
Then for words $w$ with $|w| \geq k-1$,
\[
Q_k R \xi_w = Q_k L_w R \xi_1 = Q_k (a_1\xi_{w1} + a_2 \xi_{w2^2}) =0.
\]
Thus the trace is computed as
\begin{eqnarray*}
\tr( P_{\M^\perp} Q_k P_{\M^\perp}) &=& \sum_{|w| \leq k-2} (Q_k R \xi_w, R\xi_w) \\
                    &=& \sum_{|w| = k-2}(a_1\xi_{w1}, R \xi_w) + \sum_{|w| < k-
2}(R \xi_w, R\xi_w) \\
                    &=& a_1^2 n^{k-2} + \frac{n^{k-2} - 1}{n-1}.
\end{eqnarray*}
Another computation yields,
\[
\ktil (\M^\perp) = (n-1) \lim_{k\rightarrow\infty} \frac{\tr( P_{\M^\perp} Q_k
P_{\M^\perp})}{n^k} = \frac{(n-1)a_1^2}{n^2} + \frac{1}{n^2} = r.
\]
The examples constructed show that the interval $( 1 / 9, 1 / 4 ]$ is obtained in the range of
$\ktil$ with $n = 3$ and $\prk(A) = 1$. It follows that intervals of the form $( k / 9, k / 4]$ can
be obtained with $n=3$ and $\prk(A) = k$. This completes the proof since $5k > 4$ for $k\geq
1$.
\bx

These examples are not as satisfying numerically as the decaying $n$-tuples since arbitrarily
large $n$ and pure ranks must be used. However, the associated $\fL_n$-invariant subspace is
always cyclic. They also give an indication of  how the connection with
dilation theory can be used to derive information on the ranges of the invariants.

\baselineskip=12pt

\bigskip

\begin{tabbing}
{\it E-mail address}:xx\= \kill
\noindent {\footnotesize\it Address}:
\>{\footnotesize\sc Department of Pure Mathematics}\\
\>{\footnotesize\sc University of Waterloo}\\
\>{\footnotesize\sc Waterloo, ON\quad N2L 3G1}\\
\>{\footnotesize\sc CANADA}\\
\\
{\it E-mail address}:xx\= \kill
\noindent {\footnotesize\it Current Address}:
\>{\footnotesize\sc Department of Mathematics}\\
\>{\footnotesize\sc University of Iowa}\\
\>{\footnotesize\sc Iowa City, IA\quad 52242}\\
\>{\footnotesize\sc USA}\\
\\
{\footnotesize\it E-mail address}:
\>{\footnotesize\sf dkribs@math.uiowa.edu}

\end{tabbing}

\bigskip

\bigskip

\thanks{2000 {\it Mathematics Subject Classification.} 47A13,47A20.}

\end{document}